\documentclass[12pt]{article}
\usepackage{amsmath}
\usepackage{amssymb}

\usepackage{graphicx}

\def\qed{\hfill$\Box$\par\medskip\par\relax}

\newcounter{exemplo}
\newenvironment{example}{\par\medskip\refstepcounter{exemplo}%
\noindent{\bf Example~\arabic{exemplo}.}}%
{\hfill\qed\par\medskip}

\def\eps{\varepsilon}
\def\Z{{\mathbb Z}}

\def\S{{\mathbb S}}
\def\V{{\mathcal V}}
\def\K{{\mathcal K}}
\def\M{{\mathcal M}}
\def\G{{\mathcal G}}
\def\R{{\mathbb R}}

\def\FF{{\mathfrak F}}
\def\Q{{\mathbb Q}}
\let\phi=\varphi
\def\A{{\mathfrak A}}

\def\tV{{\tilde{\mathcal V}}}
\def\tQ{{\tilde Q}}
\def\tT{{\tilde T}}

\newcommand{\eqlaw}{\stackrel{\text{\tiny law}}{=}}

\def\o{\text{\boldmath ${\omega}$}}

\def\EA{{\mathbf E}}
\def\PA{{\mathbf P}}
\def\IE{{\mathbb E}}
\def\IP{{\mathbb P}}
\def\Po{{\mathtt P}_{\omega}}
\def\Eo{{\mathtt E}_{\omega}}
\def\PoA#1{{\mathtt P}_{\omega|#1}}
\def\EoA#1{{\mathtt E}_{\omega|#1}}
\def\p{{\mathbf p}}
\def\de{{\delta}}

\def\1#1{{\bf 1}_{\{#1\}}}

\def\BRW{{branching random walk}}

\allowdisplaybreaks

\newtheorem{theo}{Theorem}[section]
\newtheorem{lmm}{Lemma}[section]
\newtheorem{df}{Definition}[section]
\newtheorem{prop}{Proposition}[section]

\def\supp{\mathop{\rm supp}}

\title{On multidimensional branching
      random walks in random environment}

\author{Francis~Comets\thanks{Partially
   supported by CNRS (UMR 7599
``Probabilit{\'e}s et Mod{\`e}les
Al{\'e}atoires''), and
by the ``R\'eseau
Math\'ematique France-Br\'esil"}$^{~,1}$ \and
 Serguei~Popov\thanks{Partially supported by CNPq (302981/02--0), CNPq/PADCT,
 and by the ``Rede Mate\-m\'atica Brasil-Fran\c{c}a"}$^{~,2}$}

\begin{document}

\maketitle

{\footnotesize
\noindent $^{~1}$Universit{\'e} Paris 7, UFR de Math{\'e}matiques,
case 7012, 2, place Jussieu, F--75251 Paris Cedex 05, France

\noindent e-mail: \texttt{comets@math.jussieu.fr}

\noindent $^{~2}$Instituto de Matem{\'a}tica e Estat{\'\i}stica,
Universidade de S{\~a}o Paulo, rua do Mat{\~a}o 1010, CEP 05508--090,
S{\~a}o Paulo SP, Brasil

\noindent e-mail: \texttt{popov@ime.usp.br}

}

\begin{abstract}
We study \BRW s in random i.i.d.\ environment in $\Z^d, d \geq 1$.
For this model, the population size cannot decrease,
and a natural definition of recurrence is
introduced. We prove a dichotomy for recurrence/transience,
 depending only on the support of the environmental law.
We give sufficient conditions for recurrence and for transience.
In the recurrent case, we study the asymptotics of the
tail of the distribution of the hitting times and
prove a shape theorem for the set of lattice
sites which are visited up to a large time. 
\\[.3cm]{\bf Short Title:}  Branching random walks in random environment
\\[.3cm]{\bf Keywords:} shape theorem, recurrence, transience, 
subadditive ergodic theorem, nestling, hitting time
\\[.3cm]{\bf AMS 2000 subject classifications:} Primary 60K37;
secondary 60J80, 82D30
\end{abstract}

\section{Introduction and results}
\label{intro}

Branching random walks in random environment
provide microscopic models for reaction-diffusion-convection
phenomena in a space inhomogeneous medium.
On the other hand, many progresses have been achieved in the last decade in the
understanding of  random walks in random environment on
$\Z^d$, see~\cite{Z}
for a recent review.
It is natural to investigate such branching walks,
and to relate the results to the non-branching case.
In this paper we continue the line of research
of~\cite{CMP,MP1,MP2}:
each particle gives birth to at least one descendant, according to
branching and jump probabilities which depend on
the location, and are given by an independent identically distributed
random field (environment). We stress here that the branching
and the transition mechanisms are {\it not\/} supposed to be
independent, and moreover, differently from~\cite{CMP,MP1,MP2},
we do not suppose that the {\it immediate\/} offsprings of a
particle jump independently. 
For an appropriate notion of recurrence and transience we prove that
either the \BRW\ is recurrent for almost every
environment or the \BRW\ is transient for almost every
environment. In addition, we show that details of the distribution of the
environment does not matter, but only its support. Though we could not
give a complete (explicit) classification in the spirit of \cite{CMP,MP1,MP2},
this is quite interesting in view of the
difficulty of the corresponding question for random walks in random
environment. For non reversible random walks (without branching) in random
environment, only few explicit results concerning recurrence/transience
are known. In the paper~\cite{La} the case of balanced environment
(i.e., the local drift is always~$0$) was studied, and it was proved
that the process is recurrent for $d\leq 2$ and transient for $d\geq 3$.
For general random walks in random environment there are conditions
sufficient for the random walk to be ballistic (and, consequently, transient),
see~\cite{Szn_PTRF02,Szn_AP03,Z}. These conditions, though, normally are not 
easily verifiable. On the other hand, for the model of the present paper,
we give explicit (and easily verifiable) conditions
for recurrence and transience, that, while failing to produce
a complete classification, nevertheless work well in many concrete
examples.
 
Also, we give a shape theorem for the set of visited
sites in the recurrent case. In terms of random walks in random environment,
this case corresponds to
nestling walks as well as to non nestling ones with
strong enough branching, and the result is once again
interesting in view of the law of large numbers for random walks,
which relies on some renewal structure and requires specific assumptions
in the current literature.

Some interesting problems, closely related to shape theorems,
arise when studying properties of tails of the distributions of
first hitting times. Here we show that \BRW{}s in random environment
exhibit very different behaviors in dimensions~$d=1$ and $d\geq 2$
from the point of view of hitting time distribution: in the
recurrent case, the annealed expectation of hitting times is always 
finite in $d\geq 2$, but it is not the case for the one-dimensional model.
Hitting times for random walks without branching in random environment
have motivated a number of papers. Tails of hitting times
distributions have a variety of different behaviors for random 
walks in random environment in dimension $d=1$, they have been
extensively studied both under the annealed law \cite{DPZ,PPZ}
and the quenched law \cite{GZ,PP}. Also, in higher dimensions there
are results on hitting times of hyperplanes, cf.~\cite{Szn_PTRF99,Szn_JEMS00}.

Many other interesting topics are left untouched
in this paper. They include shape theorems for the
transient case, questions related to the (global and local) size of
the population, hydrodynamical equations,
etc. Also, it is a challenging problem to
find the right order of decay for the tails 
of hitting times in dimension $d\geq 2$, in the recurrent case.
Finally, since in our model each particle has at least
one descendant, we do not deal at all with extinction,
which seems to be a difficult issue in random environment.

An important ingredient in our paper is the notion of seeds,
i.e., local configurations of the environment. Some seeds can create
an infinite number of particles without help from outside,
potentially enforcing recurrence. So, as opposed to random walks
without branching, the model of the present paper is in some sense
more sensible to the local changes in the environment. Together
with the fact that more particles means more averaging, this explains
why the analysis is apparently easier for the random walks with the
presence of branching.

We briefly discuss different, but related, models.
A multidimensional ($d\geq 3$) \BRW\ for which the transition probabilities
are those of the simple random walk, and the particles can branch
only in some special sites (randomly placed, with a decreasing density)
was considered in~\cite{HMP}, and several sufficient conditions
for recurrence and transience were obtained.
Dimension $d=1$ leads to more explicit results, thanks to the order structure
(see, e.g.,~\cite{CMP}). In the case $d=1$ with nearest neighbor jumps,
particles have to visit all intermediate locations, and this
fact allows to obtain some useful variational formulas \cite{GdH, BCGdH}.
The case where particles move on the tree
has a similar flavour~\cite{MP2}. The case of inhomogeneous jumps with
constant branching rate can be formulated as a tree-indexed random walk.
In this case, a complete classification of recurrence/transience is
obtained in~\cite{GM}, involving the branching rate and the spectral radius
of the transition operator. The occurrence of shape theorems in the
\BRW\ literature goes back at least to~\cite{B}.
\medskip

We now describe the model. Fix a finite set $\A\subset \Z^d$ such
that $\pm e_i \in \A$ for all $i=1,\ldots,d$, where $e_i$-s are the
coordinate vectors of $\Z^d$.
Define (with $\Z_+=\{0,1,2,\ldots\}$)
\[
 \V = \Big\{v=(v_x, x\in \A) : v_x \in \Z_+, \sum_{x\in \A}v_x \geq 1 \Big\},
\]
and for $v\in\V$ put $|v|=\sum_{x\in \A}v_x$; 
note that $|v|\geq 1$ for all $v\in\V$.
Furthermore, define $\M$ to be the set of all probability
measures~$\omega$ on $\V$:
\[
 \M = \Big\{\omega = (\omega(v), v\in\V) : \omega(v)\geq 0 \mbox{ for all }v\in\V,
                     \sum_{v\in\V}\omega(v)=1\Big\}.
\]
Finally, let $Q$ be a probability measure on $\M$. Now, for each $x\in\Z^d$ we
choose a random element $\omega_x\in\M$ according 
to the measure $Q$, independently.
The collection $\o = (\omega_x, x\in\Z^d)$ is called {\it the environment}.
Given the environment~$\o$, the evolution of the process is described in the
following way: start with one particle at some fixed
site of~$\Z^d$. At each integer time
the particles branch independently using the following mechanism:
for a particle at site $x\in\Z^d$, a random element $v=(v_y, y\in \A)$
is chosen with probability $\omega_x(v)$, and then the particle
is substituted by $v_y$ particles in $x+y$ for all $y\in \A$.

Note that this notion of \BRW\ is more general than that of \cite{CMP,MP1,MP2},
since here we do not suppose that the {\it immediate\/} descendants of
a particle jump independently (for example, we allow situations similar
to the following one [d=1]: when a particle in~$x$ generates three offsprings,
then with probability 1 two of them go to the right and the third one to the left).

We denote by $\IP, \IE$ the probability and expectation with respect to
$\o$ (in fact, $\IP=\bigotimes_{x\in\Z^d}Q_x$, where $Q_x$ are copies of $Q$),
and by  $\Po^x, \Eo^x$  the (so-called ``quenched")
probability and expectation for the process starting from~$x$ in the fixed
environment~$\o$. We use the notation
$\PA^x[\,\cdot\,] = \IE\,\Po^x[\,\cdot\,]$
for the annealed law of the branching random walk in random environment, and
$\EA^x$ for the corresponding expectation. Also, sometimes we use the symbols
$\Po,\Eo,\PA,\EA$ without the corresponding superscripts when it can create
no confusion (e.g.\ when the starting point of the process is
indicated elsewhere).

Throughout this paper, and without recalling it explicitly,
we suppose that the two conditions below are fulfilled:

\medskip
\noindent
{\bf Condition B.}
\[
Q\{\omega: \mbox{ there exists } v\in\V \mbox{ such that }
                         \omega(v)>0 \mbox{ and } |v|\geq 2\}>0.
\]

\medskip
\noindent
{\bf Condition E.}
\[
 Q\Big\{\omega: \sum_{v:v_e\geq 1} \omega(v) >0 \mbox{ for any }
 e\in \{\pm e_1,\ldots,\pm e_d\} \Big\} = 1.
\]

Condition B ensures that the model cannot be reduced to random walk without
 branching, and Condition~E
is a natural
ellipticity condition which ensures that the walk is really
$d$-dimensional. In this paper, ``elliptic'' will mean ``strictly elliptic''.
We will sometimes use the stronger
uniform ellipticity condition:

\medskip
\noindent
{\bf Condition UE.} For some $\eps_0 >0$,
\[
 Q\Big\{\omega: \sum_{v:v_e\geq 1} \omega(v) \geq \eps_0 \mbox{ for any }
 e\in \{\pm e_1,\ldots,\pm e_d\} \Big\} = 1.
\]

Due to Condition B, the population size tends to infinity, and the \BRW\
is always transient as a process on $\Z_+^{\Z^d}$. So, we
intruduce more appropriate notions of recurrence and transience.

\begin{df}
\label{def_rec}
For the particular realization of the random environment~$\o$,
the branching random walk is called recurrent if
\[
\Po^0[\mbox{the origin is visited infinitely often}]=1.
\]
Otherwise, the branching random walk is called
transient.
\end{df}
By the Markov property, the recurrence is equivalent to
\[
\Po^0[\text{the origin is visited at least once}]=1.
\]
In principle,
the above definition could depend on the starting point of the process
and on the environment~$\o$.
However, a natural dichotomy takes place:

\begin{prop}
\label{rec/tr}
We have either:
\begin{itemize}
\item[(i)] For $\IP$-almost all~$\o$, the branching random walk is recurrent, in
which case $\Po^x[\mbox{the origin is visited infinitely often}]=1$ 
for all~$x\in\Z^d$,
or:
\item[(ii)] For $\IP$-almost all~$\o$, the branching random walk is transient, in
which case $\Po^x[\mbox{the origin is visited infinitely often}]<1$ 
for all~$x\in\Z^d$.
\end{itemize}
\end{prop}

The next proposition refines the item~(ii) of Proposition~\ref{rec/tr}:
\begin{prop}
\label{trans_inf}
Let us assume that the branching random walk is transient. 
Then for $\IP$-almost all~$\o$
we have $\Po^x[\mbox{the origin is visited infinitely often}]=0$ 
for all~$x\in\Z^d$.
\end{prop}
It is plain to construct (see e.g.\ the example after the proof of Theorem~4.3
in~\cite{CMP}) environments $\o$ such that
$\Po^x[0 \mbox{ is visited infinitely often}]$ is strictly between
0 and 1.
The next example show that randomness of the  environment is essential for
our statements (and also shows, incidentally, that there is no hope to
prove Proposition~\ref{trans_inf} by arguments of the type 
``recurrence should not be
sensitive to changes of the environment in finite regions'').
\begin{example}
 Let $d=1, \A=\{-1,1\}$, and consider two measures
$\omega^{(1)}, \omega^{(2)}$:
\begin{itemize}
\item[(i)]
under $\omega^{(1)}$, with probability $2/3$ there is only one child
which is located one step to the left and  with probability $1/3$
there is
only one child which is located one step to the right;
\item[(ii)] under $\omega^{(2)}$, with probability $1/3$ there is only one child
which is located one step to the right and  with probability $2/3$
there are two children one being located to the right and the other to
the left.
\end{itemize}

If all sites $x < 0$ have the environment $\omega^{(1)}$
(we say they are of type 1) and all sites $x \geq 0$
are of type 2, we have
$\Po^x$[0 is visited infinitely often] is 1 for $x\geq 0$ and
is less than 1 for $x<0$. Changing the site $x=0$ from type 2 to type 1
turns the \BRW \ from recurrent to transient.
This example also shows that, in general,
the recurrence does depend on the environment
locally. Moreover, it shows that $\Po^0[\mbox{the origin is
visited infinitely often}]$ may be different from 0 and 1.
We will see below that, selecting randomly the environment
in an i.i.d.\ fashion, makes this \BRW\ recurrent (for this particular
example it follows e.g.\ from Theorem~\ref{suff_rec}).
\end{example}

Now, we begin stating the main results of this paper.
In Section~\ref{intro_tr_rec}
we formulate the results concerning transience and recurrence of the process,
in Section~\ref{intro_shape} we deal with questions related to (quenched
and annealed) distributions of hitting times and
shape theorems.

\subsection{Transience and recurrence}
\label{intro_tr_rec}

It is worth keeping in mind a particular example to illustrate our results.
\begin{example}
\label{exx}
With $d=2$ and $\A=\{\pm e_1, \pm e_2\}$, consider
the following $v$'s:
$v^{(1)}=\de_{e_1}$,
(with $\de$ the Kronecker symbol),
$v^{(2)}=\de_{e_2}$, $v^{(3)}=\de_{-e_1}$, $v^{(4)}=\de_{-e_2}$,
$v^{(5)}=\de_{e_1}+2\de_{e_2}+\de_{-e_1}+\de_{-e_2}$,
and the  following $\omega^{(0)}, \omega^{(+)}, \omega^{(-)}$ defined by
($0<a <1$)
\begin{eqnarray*}
 \omega^{(0)}(v^{(1)})=\frac{3}{8},\;
\omega^{(0)}(v^{(2)})=\frac{1}{4},&
\omega^{(0)}(v^{(3)})=\displaystyle\frac{1}{8}, \;
\omega^{(0)}(v^{(4)})=\frac{1}{4},& \\
 \omega^{(+)}(v^{(1)})=a,&
\omega^{(+)}(v^{(5)})=1-a, \\
\omega^{(-)}(v^{(3)})=\frac{1}{8},&
\omega^{(-)}(v^{(5)})=\displaystyle\frac{7}{8}\;,&
\end{eqnarray*}
see Figure~\ref{f_ex2}.
\begin{figure}
\centering
\includegraphics{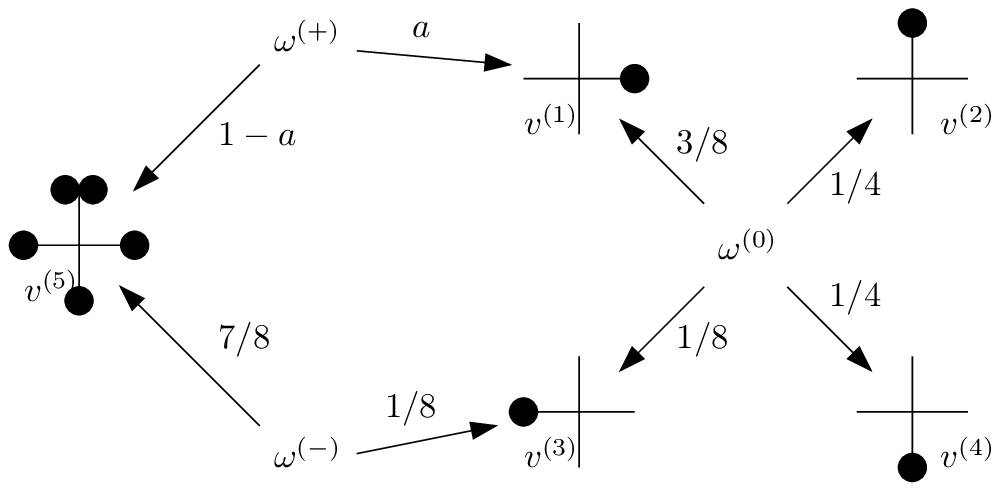}
\caption{The random environment in Example~\ref{exx}}
\label{f_ex2}
\end{figure}
Note that Conditions B and UE are satisfied.

It seems clear that the \BRW \
with $Q=Q_1$ such that
$Q_1(\omega^{(0)})=\alpha = 1-Q_1(\omega^{(-)})$
is recurrent at least for small~$\alpha$.
In fact, it is recurrent for all $\alpha \in (0,1)$.
It seems also clear that \BRW \ with $Q=Q_2$ such that
$Q_2(\omega^{(0)})=\alpha = 1-Q_2(\omega^{(+)})$ may be  recurrent or transient
depending on~$a$. In fact, for $a\leq 1/2$ the condition (\ref{sdf'}) is
fulfilled and the \BRW \ is recurrent, though Condition~L in~(\ref{eq_C_L}) 
is satisfied for $a \geq 8/9$ and  the \BRW \ is
transient (to verify~(\ref{eq_C_L}), use $s=e_1$ and $\lambda=1/3$).
But it is not so clear that the behavior
does not depend on~$\alpha$ provided that $0<\alpha<1$.
All these statements are direct applications of the
three theorems of this section.
\end{example}

Our first result states that transience/recurrence of the process only
depends on the support of the measure~$Q$, i.e., the smallest closed
subset $F \subset {\cal M}$ such that $Q(F)=1$. We recall that
$\omega$ belongs to the support iff $Q({\cal N})>0$
for all neighborhood ${\cal N}$ of $\omega$ in ${\cal M}$.
\begin{theo}
\label{trans/rec}
Suppose that the branching random walk is recurrent 
(respectively, transient) for almost all
realizations of the random environment from the distribution~$Q$.
Then for any measure~$Q'$ with $\supp Q = \supp Q'$
the process is recurrent (respectively, transient) for almost all
realizations of the random environment from the distribution~$Q'$.
(We recall that we assume that $Q'$ satisfies Condition E.)
\end{theo}

\begin{figure}
\centering
\includegraphics{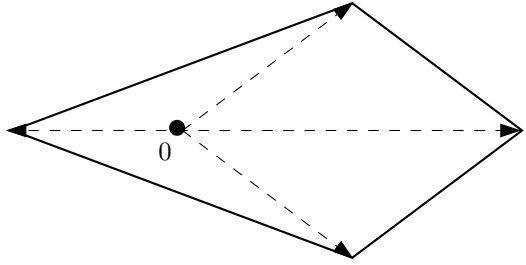}
\caption{The set $\big\{\sup_{\omega \in {\rm supp}Q}
\{\sum_v \omega(v) D(r,v)\} r ; r \in \S^{d-1}\big\}$ 
for the branching random walk (the one defined by~$Q_2$, with $a<1/2$)
from Example~\ref{exx} is the solid line}
\label{f_th1-2}
\end{figure}

Unlike the corresponding results of~\cite{CMP,MP1,MP2},
here we did not succeed in
writing down an explicit criterion
for recurrence/transience. However, sufficient condition for recurrence
or transience can be obtained (they are formulated in terms of the support
of $Q$, as they should be).
To this end, for any $v\in\V$
and any vector $r\in\R^d$, define (see Figure~\ref{f_th1-2})
\begin{equation}
\label{def_D}
D(r,v) = \max_{x\in \A : v_x \geq 1} r\cdot x,
\end{equation}
where $a\cdot b$ is the scalar product of $a,b\in\R^d$.
Let also $\|\cdot\|$ be the Euclidean norm and
$\S^{d-1} = \{a\in\R^d : \|a\|=1\}$ be the unit sphere in~$\R^d$.

\begin{theo}
\label{suff_rec}
If
\begin{equation}
\label{sdf}
\sup_{\omega\in\supp Q} \sum_{v\in\V} \omega(v) D(r,v) > 0
\end{equation}
for all $r\in\S^{d-1}$, then the branching random walk is recurrent.
Moreover, if
\begin{equation}
\label{sdf'}
\sup_{\omega\in\supp Q} \sum_{v\in\V} \omega(v) D(r,v) \geq 0
\end{equation}
for all $r\in\S^{d-1}$ and Condition~UE holds,
then the branching random walk is recurrent.
\end{theo}

Note that~(\ref{sdf'}) cannot guarantee
the recurrence without Condition~UE. To see this, consider the
following
\begin{example}
Let $d=1$, and~$Q$ puts positive weights on $\omega^{(n)}$, $n>5$,
where $\omega^{(n)}$ is described in the following way. A particle
is substituted by in mean $\frac{n}{n-1}$ offsprings
(for definiteness, let us say that it is substituted by~$2$ offsprings
with probability $\frac{1}{n-1}$ and by~$1$ offspring
with probability $\frac{n-2}{n-1}$); 
each one of the offsprings
goes to the left with probability~$1/n$, to the right with probability~$4/n$,
and stays on its place with probability $1-\frac{5}{n}$, independently.
In this case~(\ref{sdf'}) holds, but we do not have Condition~UE.
Applying Theorem~\ref{suff_tr} below (use $\lambda=1/2$), one can
see that this \BRW\ is transient.
\end{example}

\noindent
{\bf Remarks.}

(i)
Two rather trivial sufficient conditions for recurrence are: there
is $\omega\in\supp Q$ such that
\begin{equation}
\label{triv1}
\sum_{v\in\V} \omega(v) |v| = +\infty,
\end{equation}
or such that
\begin{equation}
\label{triv2}
\sum_{v\in\V} \omega(v) v_0 > 1.
\end{equation}

The proof is given after the proof of Proposition \ref{trans_inf}.

(ii)
A particular case of the model considered here is the
usual construction of the branching random walk, e.g.~\cite{CMP,MP1,MP2}:
in each~$x$,  specify the transition probabilities ${\hat P}_y^{(x)}$,
$y\in \A$, and branching probabilities ${\hat r}_i^{(x)}$, $i=1,2,3,\ldots$.
A particle in~$x$ is first substituted by~$i$ particles with
probability ${\hat r}_i^{(x)}$, then each of the offsprings jumps independently
to~$x+y$ with probability ${\hat P}_y^{(x)}$. The pairs $(({\hat
  r}_i^{(x)})_{i\geq 1}, ({\hat P}_y^{(x)})_{y\in \A})$ are chosen
according to some i.i.d.\ field on $\Z^d$.
In our notations, $\omega_x$ is a mixture of multinomial
distributions on~$\A$:
\[
   \omega_x(v) = \sum_{i\geq 1} {\hat r}_i^{(x)} 
        \mathop{\text{\sf Mult}}(i;{\hat P}_y^{(x)}, x\in\A).
\]

Note that, in this case, $D$ defined in (\ref{def_D}) is trivially
related to the local drift of the walk by
\[
 \sum_{v\in\V} \omega_x(v) D(r,v) \geq r \cdot \sum_{y \in \A}
y {\hat P}_y^{(x)}.
\]
The family of transition probabilities  ${\hat P}_y^{(x)}$,
$y\in \A$, defines a random walk in i.i.d.\ random environment
on $\Z^d$. The following definitions are essential
in the theory of such walks~\cite{Z,Zer}, they are formulated
here in the spirit of~(\ref{def_D}). 
With $\hat Q$ the common law
of $({\hat P}_y^{(x)})_{y\in \A})$, the random walk is
\begin{itemize}
 \item {\it nestling}, if for all $r \in \S^{d-1}$,
\[
\sup_{\omega\in\supp {\hat Q}}   r\cdot  \sum_{y \in \A}
y {\hat P}_y > 0;
\]
 \item {\it non nestling}, if there exists $r \in \S^{d-1}$ such that
\[
\sup_{\omega\in\supp {\hat Q}}   r\cdot  \sum_{y \in \A}
y {\hat P}_y < 0;
\]
 \item {\it marginally nestling}, if
\[
  \min_{r \in \S^{d-1}}
\sup_{\omega\in\supp {\hat Q}}   r\cdot  \sum_{y \in \A}
y {\hat P}_y = 0.
\]
\end{itemize}
Suppose now that the
 random walk in random environment
 is  nestling (or, either nestling or  marginally nestling
with Condition UE). Then, under
Condition~B, Theorem~\ref{suff_rec} implies
that the branching random walk
is recurrent, regardless of the amount of branching that
is present in the model and even though the effective drift of the
random walk can be arbitrarily large.
This extends the observation made in this case in dimension $d=1$,
Example 1 in Section 4 of~\cite{CMP}, to arbitrary dimension and
more general \BRW s. The heuristic scenario to produce such effects
remains the same: due to the nestling assumption,
 the medium develops untypical
configurations which traps the particles at some distance from the origin;
there, the branching generates an exponential number of particles,
which will balance the small probability for returning to the origin.
Indeed, the quenched large deviations rate function vanishes at 0
in the nestling case, cf.~\cite{Zer,Var}.
\qed
\medskip

Now, we turn our attention to the conditions for transience.
Define for $\omega\in\M$, $y\in \A$
\[
 \mu_y^\omega = \sum_{v\in\V} v_y \omega(v),
\]
i.e., $\mu_y^\omega$ is the mean
number of particles sent from~$x$ to~$x+y$
when the environment at~$x$ is~$\omega$.
Consider the following

\medskip
\noindent
{\bf Condition L.} There exist $s\in\S^{d-1}, \lambda>0$ such that
for all $\omega\in\supp Q$ we have
\begin{equation}
\label{eq_C_L}
 \sum_{y\in \A} \mu_y^\omega \lambda^{y\cdot s} \leq 1.
\end{equation}
We note that, by continuity, Condition~L is satisfied if and only if
(\ref{eq_C_L}) holds for $Q$-a.e.\ $\omega$.
\begin{theo}
\label{suff_tr}
Condition~L is sufficient for the transience of the branching random
walk in random
environment. Moreover, for $\IP$-a.e. $\omega$, with positive
$\Po^x$-probability the branching random walk
will not visit the half-space $\{y\in\Z^d: y\cdot s_0\leq 0\}$
-- provided that its
starting point $x$ is outside that half-space --, where
\[
 s_0 = \left\{\begin{array}{ll}
              s, & \mbox{ if } \lambda<1,\\
              -s, & \mbox{ if } \lambda >1 .
              \end{array}
          \right.
\]
(As we will see below, Condition~L cannot be satisfied
with~$\lambda=1$.) Furthermore, the number of visits of the \BRW\
to the half-space is a.s.\ finite.
\end{theo}

In~\cite{CMP,MP1,MP2} it was shown that, if~$L_0=1$, 
conditions analogous to Condition~L
are sufficient and {\it necessary\/} for transience in cases when the
branching random walk in random environment lives on the one-dimensional lattice or
on a tree (in particular, by repeating the argument of~\cite{CMP}, it is
not difficult to prove that for the present model with~$L_0=1$ in dimension~$1$,
Condition~L is necessary and sufficient for transience). 
On $d$-dimensional lattice ($d\geq 2$) or even in
dimension~$1$ with $L_0>1$ the question whether
Condition~L is necessary for transience remains open.

\subsection{Hitting times and asymptotic shape in the recurrent case}
\label{intro_shape}
For the process starting from one particle at~$x$, let us denote by
$\eta_n^{x}(y)$  the number of particles in~$y$ at time~$n$, and
by $B_n^x$  the set of all sites visited by the process
before time~$n$. Also, denote by $T(x,y)$ the moment of hitting~$y\neq x$.
For the formal definition of those quantities,
see Section~\ref{sec:constr}, although here
we do not need to construct simultaneously
all the \BRW s from all the possible starting points $x$.

First, we are going to take a closer look at the hitting times
$T(0,x)$ for recurrent branching random walks. It immediate to
note that the recurrence is equivalent to
$\PA[T(0,x)<\infty \mbox{ for all }x]=1$. So, for the recurrent
case it is natural to ask how fast the recurrence occurs,
i.e., how fast (quenched and annealed) tails of the
distribution of $T(0,x)$ decrease. For the (quenched) asymptotics
of $\Po[T(0,1)>n]$ in dimension~$1$, we have the following
result:
\begin{prop}
\label{quench_dim1}
Suppose that $d=1$ and the branching random walk 
in random environment is recurrent.
Then, for $\IP$-almost all environments there exist $n^*=n^*(\o)$
and~$\kappa>0$ such that
\begin{equation}
\label{eq_quench_dim1}
  \Po[T(0,1)>n] \leq e^{-n^{\kappa}}
\end{equation}
for all $n\geq n^*$.
\end{prop}
This result follows from a more general fact that will be proved
in the course of the proof of Theorem~\ref{shape_t}, case $d=1$
(see the remark just below the formula~(\ref{oc_strexp})). 
Moreover, the following example
shows that, for the class of recurrent one-dimensional
branching random walks in random environment,
the right order of decay of $\Po[T(0,1)>n]$ is indeed
stretched exponential.

\begin{example}
\label{example_qdecay}
We consider $d=1$, $\A=\{-1,1\}$, and
suppose that~$Q$ gives weights $1/3$ to the points $\omega^{(1)},
\omega^{(2)},\omega^{(3)}$, which are described as follows.
Fix a positive $p<1/82$; there is no branching in $\omega^{(1)},
\omega^{(2)}$, and $\omega^{(1)}$ (respectively, $\omega^{(2)}$),
sends the particle to the left (respectively, to the right) 
with probability~$p$ and to the right (respectively, to the left) 
with probability~$1-p$. In the sites with $\omega^{(3)}$, the
particle is substituted by~$1$ (respectively, $2$) offsprings with
probability~$2p$ (respectively, $1-2p$); those then jump
independently to the right or to the left with equal probabilities.
By Theorem~\ref{suff_rec}, this branching random walk is
recurrent.

Abbreviate for a moment $a=\ln\frac{1-p}{p}$; note that
\[
\IP\big[\omega_{-x}=\omega^{(1)}, \omega_x=\omega^{(2)}
\mbox{ for }x\in (0,a^{-1}\ln n], \omega_0=\omega^{(1)}\big] =
  n^{-2a^{-1}\ln 3}.
\]
Clearly, $p<1/82$ implies that $2a^{-1}\ln 3<1/2$. This means that
 a typical environment~$\o$ will contain a translation of the
trap considered above in the box $[-n^{1/2},0]$, i.e., there
is an interval $[b-a^{-1}\ln n, b+a^{-1}\ln n]\subset [-n^{1/2},0]$
such that $\omega_x=\omega^{(1)}$ for $x\in [b-a^{-1}\ln n, b]$
and $\omega_x=\omega^{(2)}$ for $x\in (b, b+a^{-1}\ln n]$.
For such an environment, first, the initial particle goes
straight to the trap (without creating any further particles on its way) 
with probability at least $p^{n^{1/2}}$,
and then stays there with a probability bounded away from~$0$
(note that the depth of the trap is $\ln n$, and this is enough to
keep the particle there until time~$n$ with a good probability). 
This shows that $\Po[T(0,1)>n]\geq C e^{-n^{1/2}\ln p^{-1}}$.
\end{example}

One can construct other one-dimensional examples of this type;
the important features are:
\begin{itemize}
\item[(i)] there are $\omega$-s from $\supp Q$ without branching and with drifts
pointing to both directions, so that traps are present;
\item[(ii)] all $\omega$-s from $\supp Q$ have the following
property: with a positive probability the particle does not
branch, i.e., it is substituted by only one offspring; this
permits to a particle to cross a region without necessarily
creating new ones.
\end{itemize}

In dimensions $d\geq 2$, finding the right order of decay of $\Po[T(0,1)>n]$
is, in our opinion, a challenging problem. For
now, we can only conjecture that it should be exponential
(observe that the direct attempt to
generalize the above example to $d\geq 2$ fails, since
creating a trap with a logarithmic depth has higher cost).
As a general fact, the annealed bound of Theorem~\ref{t_upbound_ht}
below is also a quenched one. This is the only rigorous result
concerning the quenched asymptotics of $\Po[T(0,1)>n]$ we can 
state in the case $d\geq 2$; 
we believe, however, that it is far from being precise.

Now, we turn our attention to the annealed distribution of hitting times.
\begin{theo}
\label{t_upbound_ht}
Let $d\geq 1$.
For any $x_0\in\Z^d$ there exists $\theta=\theta(x_0, Q)$ such that
\begin{equation}
\label{eq_hit_time}
\PA[T(0,x_0)>n] \leq \exp\{-\theta\ln^d n\}
\end{equation}
for all~$n$ sufficiently large.
\end{theo}

Define $\G\subset\M$ to be the set of $\omega$-s without
branching, i.e.,
\[
 \G = \Big\{\omega\in\M : \sum_{v\in\V: |v|=1} \omega(v)=1\Big\}.
\]
In other words, if at a given~$x$ the environment belongs
to~$\G$, then the particles in~$x$ only jump, without creating
new particles. Also, for any $\omega\in\G$, define the drift
\[
 \Delta_\omega = \sum_{x\in\A} x\omega(\delta_x).
\]
The following result shows that Theorem~\ref{t_upbound_ht}
gives in some sense the best possible bounds for the 
tail of the hitting time distribution that are valid for the class
of recurrent branching random walk in random environment.

\begin{theo}
\label{t_lobound_ht}
Suppose that $Q(\G)>0$ and that the origin belongs to the 
interior of the convex hull of $\{\Delta_\omega :
\omega\in \G\cap \supp Q\}$. Then for 
any $x_0\in\Z^d$ there exists $\theta'=\theta'(x_0, Q)$ such that
\begin{equation}
\label{lowerbound_ht}
\PA[T(0,x_0)>n] \geq \exp\{-\theta'\ln^d n\}
\end{equation}
for all~$n$ sufficiently large.
\end{theo}

  From Theorems~\ref{t_upbound_ht} and~\ref{t_lobound_ht}
there is only a small distance to the following remarkable
fact: the implication
\[
 \big(\PA[T(0,x)<\infty \mbox{ for all }x]=1\big) \implies 
   \big(\EA T(0,x)<\infty \mbox{ for all }x\big)
\]
is true for $d\geq 2$ but is false for $d=1$. To see this, it is
enough to know that the constant~$\theta'$ in~(\ref{lowerbound_ht})
can be less than~$1$ in dimension one. Consider the
following example:
\begin{example}
\label{ex_infexp}
Once again, we suppose that $\A=\{-1,1\}$ and
$\supp Q$ consists of three points $\omega^{(1)},
\omega^{(2)},\omega^{(3)}$, with $Q(\omega^{(1)})=\alpha_1$,
 $Q(\omega^{(2)})=\alpha_2$,  $Q(\omega^{(3)})=1-\alpha_1-\alpha_2$.
We keep the same $\omega^{(1)},\omega^{(2)}$ from Example~\ref{example_qdecay},
and let $\omega^{(3)}(\delta_1+\delta_{-1})=1$. It is immediate
to obtain from e.g.\ Theorem~\ref{suff_rec} that this branching random walk 
in random environment is recurrent. Abbreviate $a=\ln \frac{1-p}{p}$
and let 
\begin{eqnarray*}
H &=& \{\omega_x=\omega^{(1)}\mbox{ for }x\in [-2a^{-1}\ln n, -a^{-1}\ln n],\\
&& ~~~~~~~~~~~~~~~~~~~~~~~
\omega_x=\omega^{(2)}\mbox{ for }x\in (-a^{-1}\ln n,0]\};
\end{eqnarray*}
then $\IP[H]=(\alpha_1\alpha_2)^{a^{-1}\ln n}$. Now, on~$H$
there is a trap of depth~$\ln n$ just to the left of the origin,
so for such environments the quenched expectation of~$T(0,1)$ over 
all paths which visit the site $(-a^{-1}\ln n)$ before the site~$1$
is at least~$Cn$, so we have
\[
\EA T(0,1) \geq \int_{H}\Eo T(0,1) d\IP \geq C n \IP[H]
   = Cn^{1-a^{-1}\ln (\alpha_1\alpha_2)^{-1}} \to \infty
\]
when $a^{-1}\ln (\alpha_1\alpha_2)^{-1}<1$ (or equivalently,
$p<(1+(\alpha_1\alpha_2)^{-1})^{-1}$). Here we could use also the
same branching random walk of Example~\ref{example_qdecay}
(with $p<1/10$); note however that in the present example
we could allow sites where particles always branch.
\end{example}

Now, we pass to a subject closely related to hitting times,
namely, we will study the set of the sites visited by time~$n$
(together with some related questions).
Recall that
\[
 B_n^x = \{y\in\Z^d : \mbox{ there exists } m\leq n \mbox{ such that }
             \eta_m^x(y)\geq 1\}.
\]
Also, we define ${\bar B}^x_n$ as the set of sites that contain at least one
particle at time~$n$, and ${\tilde B}^x_n$ as the set of sites
that contain at least one particle at time~$n$ and will always do
in future:
\begin{eqnarray*}
{\bar B}^x_n & = & \{y\in\Z^d : \eta_n^x(y)\geq 1\}, \\
{\tilde B}^x_n & = & \{y\in\Z^d : \eta_m^x(y)\geq 1\mbox{ for all }m\geq n\}.
\end{eqnarray*}
Evidently, ${\tilde B}^x_n \subset {\bar B}^x_n  \subset B^x_n$ 
for all~$x$ and~$n$.

When dealing with the shape results for ${\tilde B}^x_n$ and
${\bar B}^x_n$, we will need the following ``aperiodicity'' condition,
where we use the standard notation
$\|x\|_1=|x^{(1)}|+\cdots+|x^{(d)}|$
for $x=(x^{(1)},\ldots,x^{(d)})\in\Z^d$:

\medskip
\noindent
{\bf Condition A.} There exist $ x \in \A , v \in \V$ with $\|x\|_1$
even and $v_x \geq 1$
such that $Q \{ \omega \in \M : \omega(v) >0 \}>0$.

\medskip

For any set $M\subset\Z^d$ we define the set $\FF(M)$ by ``filling the spaces''
between the points of~$M$:
\[
 \FF(M) = \{ y + ({-1/2}, {1/2}]^d: y \in M\}\subset\R^d.
\]

We only deal with the recurrent case here, leaving the more delicate,
transient case for
further research.

\begin{theo}
\label{shape_t}
Suppose that the branching random walk in random environment is recurrent
and Condition~UE holds.

Then there exists a deterministic compact convex set $B\subset \R^d$
with 0 belonging to the interior of $B$,
such that $\PA$-a.s. (i.e., for $\IP$-almost all~$\o$ and  $\Po$-a.s.),
we have for any $0<\eps<1$
\[
(1-\eps)B \subset
 \frac{\FF(B_n^0)}{n}
\subset (1+\eps) B
\]
for all~$n$ large enough.

If, in addition, Condition~A holds, then
the same result -- with the same limiting shape~$B$ --
holds for ${\bar B}^0_n$ and ${\tilde B}^0_n$.
\end{theo}

It is straightforward to note that~$B$ is a subset of the convex hull
of~$\A$; also, since $B_n^x\eqlaw B_n^0+x$, the same result
holds for the process starting from~$x$.
 As opposed to the results of the previous section,
the limiting shape~$B$ does not only 
depend on the support of~$Q$, see the example below:
\begin{example}
Let $d=1$ and $\A=\{-2,-1,0,1,2\}$. 
Then, put $v^{(1)}=\delta_{-1}+\delta_{0}+\delta_1$,
$v^{(2)}=\delta_{-2}+\delta_{-1}+\delta_{0}+\delta_1+\delta_2$, 
$\omega^{1}=\delta_{v^{(1)}}$,
$\omega^{2}=\delta_{v^{(2)}}$, and $Q(\omega^{1})=1-Q(\omega^{2})=\alpha$. Then,
it is quite elementary to obtain that $B=[-(2-\alpha),2-\alpha]$, 
i.e., the asymptotic shape
depends on~$Q$ itself, and not only on the support of~$Q$.
\end{example}

Another interesting point about Theorem~\ref{shape_t} is that
the shape $B$ is convex,  but one finds easily examples -- as the one
below -- where it is not strictly convex.

\begin{example}
\label{ex_flatedge}
With $d=2$ and $\A=\{\pm e_1, \pm e_2\}$, consider
$v^{(1)}=\de_{e_1}+\de_{e_2}$,
$v^{(2)}=\de_{e_1}+\de_{e_2}+\de_{-e_1}$,
$v^{(3)}=\de_{e_1}+\de_{e_2}+
\de_{-e_2}$,
and the  following $\omega^{(1)}, \omega^{(2)}$ defined by
\[
\begin{array}{c}
 \omega^{(1)}(v^{(1)})= \omega^{(2)}(v^{(1)})=
 \omega^{(1)}(v^{(2)})=
 \omega^{(2)}(v^{(3)})=\displaystyle\frac{2}{5},\\
\omega^{(1)}(v^{(3)})=\omega^{(2)}(v^{(2)}) =\displaystyle\frac{1}{5},
\end{array}
\]
see Figure~\ref{f_ex5}.
\begin{figure}
\centering
\includegraphics{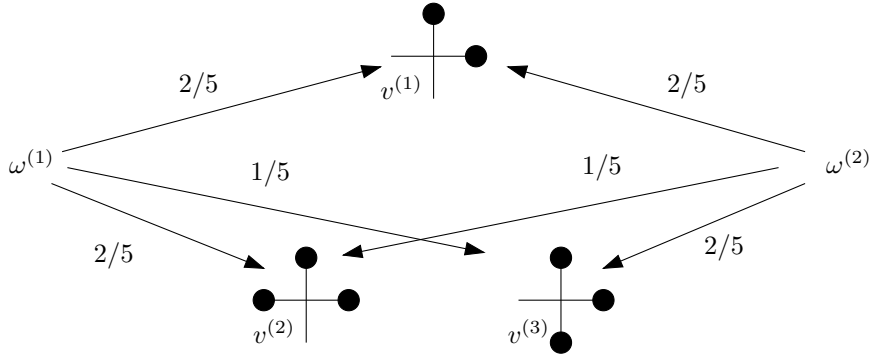}
\caption{The random environment in Example~\ref{ex_flatedge}}
\label{f_ex5}
\end{figure}
Then, with $Q(\omega^{(1)})=1-Q(\omega^{(2)})=\alpha$ (with $\alpha \in
(0,1)$), the \BRW \ is recurrent by Theorem \ref{suff_rec}.
For arbitrary $\alpha \in (0,1)$,
$B_n^0 \cap \Z_+^2=\{(x_1,x_2): x_1,x_2 \in \Z_+,
x_1+x_2\leq n\}$, and so
\[
B \cap \R_+^2=\big\{(x_1,x_2): x_1,x_2 \geq 0,
x_1+x_2\leq 1\big\},
\]
and~$B$ has a flat edge.
\end{example}

\section{Some definitions and preliminary facts}
First, let us introduce some more basic notations:
for $x=(x^{(1)},\ldots,x^{(d)})\in\Z^d$ write
\[
 \|x\|_{\infty} = \max_{i=1,\ldots,d} |x^{(i)}|.
\]
Define $L_0$ to be the maximal  jump length, i.e.,
\[
 L_0 = \max_{x\in \A} \|x\|_{\infty},
\]
and let $\K_n$ be the $d$-dimensional cube of size $2n+1$:
\[
 \K_n = [-n,n]^d = \{x\in\Z^d: \|x\|_{\infty}\leq n\}.
\]

For~$\omega\in\M$ and $V\subset\V$, sometimes we will use notations
like $\omega(v\in V)$ or even $\omega(V)$ for $\sum_{v\in V}\omega(v)$.

\subsection{Induced random walks}
\label{sec:induced}
It is most natural to connect the branching random walk in  random environment
with random walks in  random environment.
Defining
\[
\tV= \big\{ (v,\kappa): v \in \V, \kappa {\rm \ probability \ measure\ on \ }
\{y: v(y)\geq 1\} \big\},
\]
we consider some probability  measure $\tQ$ on $\tV$ with marginal $Q$
on  $\V$. An i.i.d.\ sequence
$\tilde \o = ((\omega_x, \kappa_x), x \in \Z^d)$  with distribution~$\tQ$
defines our branching random walk as above, coupled with
a random walk in random (i.i.d.) environment with transition
probability
\[
p_x(y)= \sum_{v \in \V} \omega_x(v) \kappa_x(y)
\]
from $x$ to $x+y$. In words, we pick randomly one of the children in
the branching random walk. We call this  walk, the $\tQ$-{\it induced
random walk\/} in  random  environment. Here are some natural choices
(in the examples below~$\kappa$ does not depend on~$\omega$):
\begin{itemize}
\item[(i)] uniform: $\kappa$ is uniform on the locations
$\{x \in \A: v_x\geq 1\}$;
\item[(ii)] particle-uniform: $\kappa(y)$ is proportional to the number of
particles sent by~$v$ to~$y$;
\item[(iii)] $r$-extremal, $r\in\S^{d-1}$:
 $\kappa$ is supported on the set of
$x$'s maximizing $r\cdot x$ with $x \in \A, v_x \geq 1$.
\end{itemize}

The following proposition is a direct consequence of Theorem~\ref{suff_rec}:
\begin{prop}
\label{p_nestl}
If the branching random walk in random environment is transient,
then any induced random walk is either non nestling or marginally nestling.
Moreover, if Condition~UE holds, then any induced random walk is
non nestling.
\end{prop}
Note, however, that one can easily construct examples of recurrent \BRW{}s such that
any induced random walk is non nestling, i.e., 
the converse for Proposition~\ref{p_nestl}
does not hold. For completeness we give the following
\begin{example}
Let $d=1$ and~$Q=\delta_\omega$, where~$\omega$ sends 1 particle to the left
with probability~$1/3$, and 5 particles to the right with probability~$2/3$.
Then, clearly, any induced random walk is non nestling. To see that this \BRW\
is recurrent it is enough to obtain by a simple computation that a mean number
of grandchildren that a particle sends to the same site in 2 steps is strictly
greater than 1.
\end{example}
\noindent
See also Example~2 of \cite{CMP}.

\subsection{Seeds}

In the next definition we introduce the notion of $(A,H)$-{\em seed},
which will be
frequently used throughout the paper.

\begin{df}
\label{def_seed}
Fix a finite set $A \subset \Z^d$ containing 0,
and $H_x \subset \M$ with $Q(H_x)>0$ for all $x \in A$. With
$H=(H_x, x \in A)$, the couple  $(A,H)$ is called a seed.
We say that $\o$ has a $(A,H)$-seed at $z \in \Z^d$
 if
\[
\omega_{z+x} \in H_x  \mbox{ ~for all~ } x \in A,
\]
and that $\o$ has a $(A,H)$-seed in the case $z=0$.
We call $z$ the center of the seed.
\end{df}
\begin{lmm}
\label{visit-seed}
With probability~$1$
the branching random walk visits infinitely many distinct $(A,H)$-seeds
(to visit the seed means to visit the site where the seed
is centered).
\end{lmm}

\noindent
{\it Proof.} By the ellipticity
Condition E, the uniform induced random walk is  elliptic, so for
{\it every\/} environment it cannot stay forever in a finite subset.
Take~$n$ such that $A \subset\K_n$, and partition
the lattice $\Z^d$ into translates of $\K_n$.
Since the environment is i.i.d.,
we can construct the induced random walk by choosing randomly the environment
in a translate of~$K_n$ at the first moment when the walk enters this
set. If $Q(H_x)>0$ for all $x \in A$, then by Borel-Cantelli lemma,
infinitely many
of those translates contain the desired seed, and by using Condition~E,
it is elementary to show that infinitely many seed centers will be
visited.
\qed

For a particular realization of the random environment~$\o$, we
define the branching random walk restricted on set~$M\subset\Z^d$ simply by
discarding all particles that step outside~$M$, and write~$\PoA{M}, \EoA{M}$
for corresponding probability and expectation.

\begin{df}
\label{rec_seed}
Let $U,W$ be two finite subsets of~$\Z^d$ with $0 \in W \subset U$.
Let $\p$ be a probability distribution on $\Z_+$
 with mean larger than 1, i.e.,
$\p=(p_0, p_1,p_2,\ldots)$ with $p_i\geq 0$,
$\sum p_i=1$, $\sum ip_i>1$.
An $(U,H)$-seed is called $(\p,W)$-{\em recurrent} if for any $\o$
such that $\omega_x\in H_x, x\in U$, we have
\[
\PoA{U}^y[ W {\rm \ will\ be\ visited\ by\ at\ least\ } i {\rm \
 ``free"
\ particles}]
   \geq \sum_{j=i}^\infty p_j
\]
for all $i \geq 1$ and all $y \in W$. By ``{\em free}'' particles
we mean that none is the descendant of
another one. To shorten the notations,
in the case~$W=\{0\}$ we simply say that the seed is $\p$-recurrent.
\end{df}
Note that, by definition of the restricted branching random walk,
the above probability
depends on the environment inside $U$ only.

\medskip
\noindent
The next lemma shows the relevance of $\p$-recurrent seeds.

\begin{lmm}
\label{l_rec_seed_rec}
Suppose that there exists an
$(U,H)$-seed that is $(\p,W)$-re\-cur\-rent.
Then this implies the recurrence of the \BRW\ for a.e.\ environment~$\o$.
\end{lmm}

\noindent
{\it Proof.}
We first prove this lemma when $W=\{0\}$.
It is not difficult to see that each $\p$-recurrent
seed gives rise to a super-critical Galton-Watson branching process.
More precisely,  by assumption, for all $\o$ in a set of positive
$\IP$-probability,
the \BRW\ starting at~$0$ generates
$i$ ``free'' particles with $\Po^0$-probability larger or equal
to~$p_i$,
and the particles will visit 0 before exiting $U$. Since
they are ``free'', they  are themselves walkers starting
from the site 0, each one will generate independently new
visits at 0 which number is stochastically larger than $\p$. Hence the
number of such visits at 0 dominates a super-critical Galton-Watson
branching process with offspring distribution $\p$, in the sense
that when the Galton-Watson process
survives forever, then in the original process~$0$ is visited
infinitely often.
On the other hand, by Lemma~\ref{visit-seed},
an infinite number of such seeds will be visited.
By construction, the Galton-Watson processes mentioned above
are independent for nonoverlapping seeds, so almost surely at least
one of them will generate
an infinite number of visits to the center of its seed,
thus sending also an infinite number of particles to~$0$, which proves the
recurrence for all $\o$'s under consideration.

The proof is easily extended to the case of a general~$W$.
\qed
\medskip

\noindent
{\it Proof of Proposition~\ref{rec/tr}.}

\noindent
Assume that the event $\{
\Po^0[\mbox{the origin is visited infinitely often}]=1\}$ has
positive $\IP$-probability. Then, by Condition~B, this event intersected
with $\{\mbox{there exists }x \in \Z^d: \omega_x(v: |v|\geq 2)>0\}$
has also positive $\IP$-probability. Fix $\o$ in the intersection,
and also in the support of $\IP$.
By Condition~E, the following happens with positive
$\Po^0$-probability: one particle (at
least) of the
\BRW\ reaches this branching site $x$, and is then substituted by 2 particles
(at least), each of them eventually hitting the origin. Since the
number of visits to 0 is infinite, we get by Borel-Cantelli lemma that
\[
  \Po^0[\mbox{the origin is visited by (at least)
two free particles}]=1
\]
(recall that by ``free'' particles, we mean that none is the descendant of
  the other; we do not require that they visit 0 at the same moment).
Then, we can take~$t$ large enough so that
\[
  \Po^0[\mbox{the origin is visited by (at least)
two free particles before time } t ]> 3/4.
\]
Since the jumps are bounded, this probability is equal to
\[
\PoA{{\K}_{tL_0}}^0[\mbox{the origin is visited by (at least)
two free particles before time } t ],
\]
which depends only on $\omega_x, x \in U:={\K}_{tL_0}$.
By continuity, we can choose now small neighborhoods $H_x$ of
$\omega_x, x \in U$, such that
\begin{multline*}
{\mathtt P}_{\omega'|{\K}_{tL_0}}^0
[\mbox{the origin is visited by (at least)
two free}\\\mbox{ particles before time } t ]>3/4
\end{multline*}
for all $\omega'$ with $\omega_x' \in H_x, x \in U$. By the support
condition it holds that $Q(H_x)>0$, and we see that
the $(U,H)$-seed  is $\p$-recurrent, with $\p=(1/4, 0, 3/4, 0,0, \ldots)$.
From Lemma~\ref{l_rec_seed_rec}, we conclude that the \BRW \ is
recurrent for $Q$-a.e.\ environment. Therefore the set of
 recurrent $\o$ has probability 0 or 1.
\medskip

On the other hand, it is clear by ellipticity that
\[
\o \mbox{\  recurrent\ } \iff \;
\Po^0[x \mbox{ is visited infinitely often}] = 1,
\]
for all $x \in \Z^d$. Since the set of $\o$'s such that
$\Po^x[0 \mbox{ is visited infinitely often}]=1$ is a shift of the
previous one, it has the same probability. \qed
\medskip

\noindent
{\it Proof of Proposition~\ref{trans_inf}.}
Assume that with positive $\IP$-probability,
\[
\Po^x[\mbox{the origin is visited infinitely often}]>0
\]
for some~$x\in\Z^d$, and fix such an $\o$ in the support of $\IP$.
Then, by Condition~E,  the inequality holds for all $x \in \Z^d$, and
in view of Condition~B, we can assume without loss of generality that
\begin{equation} 
\label{neige}
\omega_x(v: |v|\geq 2)>0 .
\end{equation}
By  Condition E, $x$ is visited infinitely often a.s.\ on the event
$E:=$\{0  is visited infinitely often a.s.\}. Together with~(\ref{neige})
this implies that~$x$ is visited by infinitely many free
particles a.s.\ on the event~$E$. With $\beta=\Po^x[E]>0$,
fix some integers $K, t$ such that $K\beta /2>1$ and
\[
\Po^x[\mbox{at least $K$ free particles visit $x$ before time }
  t]> \frac{\beta}{2} .
\]
We note that this probability depends only on~$\o$
inside $U:=x+\K_{tL_0}$, hence it is equal to the
$\PoA{U}^x$-probability of the event under consideration.
By continuity, we can choose small neighborhoods $H_y$ of $\omega_y, y
\in U$, such that
\[
{\mathtt P}_{\omega'|U}
[\mbox{at least $K$ free particles visit $x$ before time }
  t]>\frac{\beta}{2}
\]
for all $\o'$ with $\omega'_y \in H_y, y \in U$. We see that the
$(U,H)$-seed
is $(\p,\{x\})$-recurrent with $p_K=\beta/2=1-p_0$, and has a positive
$\IP$-probability since $\o$ is in the support of this measure.
From Lemma~\ref{l_rec_seed_rec}, we conclude that the \BRW \ is
recurrent, which completes the proof.
\qed

\medskip

We conclude this section by proving the sufficiency (for the recurrence)
of the conditions~(\ref{triv1}) and~(\ref{triv2}).
In the latter case, we define  $U=W=\{0\}$ and~$\p$ by
\[
\sum_{j=i}^\infty p_j =
\inf( \omega'\{v: v_0\geq i\}; \omega' \in {\cal N}),\; i \in \Z_+,
\]
where ${\cal N}$ is a neighborhood of $\omega$. As
 ${\cal N} \searrow \{\omega\}$,  the mean of $\p$ increases by
 continuity to $\sum_v \omega(v)|v|>1$
 in view of (\ref{triv2}).
Choosing  ${\cal N}$ small
enough so that the mean of $\p$ is larger than 1, the seed
$(\{0\}, {\cal N})$ is $\p$-recurrent.

In the former case,
we let $U=\A, W=\{0\}$, and for $x \in \A\setminus \{0\}$,
\[
H_x=
\Big\{\omega': \sum_{v:v_e\geq 1} \omega'(v) \geq \eps  \mbox{ for all }
 e\in \{\pm e_1,\ldots,\pm e_d\} \Big\},
\]
where we fix $\eps>0$ small
enough so that $Q(H_x)>0$. Fix $a > \eps^{-dL_0-1}$, and a
distribution ${\bf q}$ on $\Z_+$ which is stochastically smaller than
the distribution of $|v|$ under $\omega$ and has mean
at least $a$.
By continuity again, and  in view of (\ref{triv1}), the set
\[
H_0= \{ \omega': \mbox{distribution of $|v|$ under $\omega'$
is stochastically larger than } {\bf q} \}
\]
is a neighborhood of $\omega$. Now, if $\o$ has a  $(U,H)$-seed in 0,
a walker starting at 0 has a number~$N$
of offsprings stochastically larger than ${\bf q}$, each of which
being able to walk back to 0 simply
by ellipticity. So, we see that given $N$ the number of
free visits to 0 without exiting $U$ dominates a binomial
distribution ${\cal B}(N, \eps^{dL_0+1})$, and that, unconditionally,
it dominates the mixture $\p$ of such binomials for $N \sim {\bf q}$.
Since~$\p$ has mean larger than 1,
the seed $(U,H)$ is $\p$-recurrent, and it has positive
$\IP$-probability. This ends the proof.
\qed

\subsection{Formal construction of the process and subadditivity}
\label{sec:constr}
Recall that,
for the process starting from one particle at~$x_0$, the variable
$\eta_n^{x_0}(x)$ is the number of particles in~$x$ at time~$n$.
Given~$\o$, for all $x\in\Z^d$, consider an i.i.d.\ family
$v^{x,i}(n)$, $i=1,2,3,\ldots$, $n=0,1,2,\ldots$, of random elements of~$\V$, with
$\Po[v^{x,i}(n)=v]=\omega_x(v)$
(with a slight abuse of notations, we will still write $\Po^x$ for the
forthcoming construction for a fixed $\o$, and
$\PA^x[\,\cdot\,] = \IE\,\Po^x[\,\cdot\,]$).
Now, the idea is to construct
the {\it collection} of branching random walks {\it indexed by the
  position of the initial particle}, using the same realization of
$(v^{x,i}(n), x\in\Z^d, i=1,2,3,\ldots, n=0,1,2,\ldots)$.

To this end, consider first the process beginning at the origin, and
put
$\eta_0^0(0)=1$,
$\eta_0^0(y)=0$ for $y\neq 0$.
Inductively,  define for $y \in \Z^d$ and $n \geq 0$:
\begin{equation}
\label{coupling}
 \eta_{n+1}^0(y) = \sum_{x:y\in \A+x} \sum_{i=1}^{\eta_n^0(x)} v_{y-x}^{x,i}(n).
\end{equation}

Define $T(0,y)$ to be the first moment when a particle enters~$y$, provided
that the process started from~$0$, i.e.,
\[
 T(0,y) = \inf\{n\geq 0 : \eta_n^0(y)\geq 1\},
\]
and $T(0,y)=+\infty$ if there exists no such~$n$.
Now, the goal is to define $\eta_n^z$ for $z\neq 0$. We distinguish
two cases.

If $T(0,z)=+\infty$,
we proceed as before, i.e., put $\eta_0^z(z)=1$,
$\eta_0^z(y)=0$ for $y\neq z$, and
\begin{equation}
\label{coupling'}
 \eta_{n+1}^z(y) = \sum_{x:y\in \A+x} \sum_{i=1}^{\eta_n^z(x)} v_{y-x}^{x,i}(n).
\end{equation}
When $m_0:= T(0,z) <\infty$, we put $\eta_0^z(z)=1$,
$\eta_0^z(y)=0$ for $y\neq z$, and
\begin{equation}
\label{coupling''}
 \eta_{n+1}^z(y) = \sum_{x:y\in \A+x} \sum_{i=1}^{\eta_n^z(x)} v_{y-x}^{x,i}(n+m_0).
\end{equation}
Define also
\[
 T(z,y) = \inf\{n\geq 0 : \eta_n^z(y)\geq 1\}.
\]

Note that the set~$B_n^x$ can now be defined by
$B_n^x = \{y: T(x,y)\leq n\}$.
The following lemma will be very important in the course of the proof
of Theorem~\ref{shape_t}:
\begin{lmm}
\label{l_subadd}
For any $y,z\in\Z^d$ and for\/ {\em all} realizations of
$(v^{u,i}(n), u\in\Z^d, i=1,2,3,\ldots,n=0,1,2,\ldots)$ it holds that
\begin{equation}
\label{e_subadd}
 T(0,z) + T(z,y) \geq T(0,y).
\end{equation}
\end{lmm}

\noindent
{\it Proof.} The inequality (\ref{e_subadd}) is obvious when  $T(0,z)=+\infty$,
 so
we concentrate on the case $m_0:=T(0,z)<\infty$. In this case, by induction,
it is immediate to prove that $\eta_n^z(x) \leq \eta_{n+m_0}^0(x)$ for all
$x\in\Z^d$ and all $n\geq 0$, which, in turn, shows~(\ref{e_subadd}).
\qed

\medskip
\noindent
{\bf Remark.} For the present model we failed to construct a coupling
such that $T(x,y)+T(y,z)\geq T(x,z)$ holds for all $x,y,z\in\Z^d$.
In absence of such a coupling
we need  to use a stronger version of the Subadditive Ergodic Theorem,
a variant of Theorem~\ref{t_subadd} below.
An example of ``branching-type'' model for
which such a coupling does exist can be found in~\cite{AMP}.

\section{Proofs: Recurrence/transience}
\label{sec:pr-rec/tr}

\subsection{Proof of Theorem~\ref{trans/rec}}

We need some preparations. The following lemma complements Lemma
\ref{l_rec_seed_rec}.

\begin{lmm}
\label{l_rec_seed_exist}
Suppose that the \BRW \ in random environment from~$Q$
is recurrent.
Then there exist $\p$, $m\geq 1$,
and a collection $H=(H_z\subset\M, z\in\K_{mL_0})$
such that $Q(H_z)>0$ for all $z\in\K_{mL_0}$, and such that
the $(\K_{mL_0},H)$-seed is $\p$-recurrent.
\end{lmm}

In fact, the reader probably has noticed that a similar result was already
proved in the course of the proof of Proposition~\ref{rec/tr}. However, for later
purposes, we will construct this seed in a more explicit way 
(see Definition~\ref{good_seed} below).

\medskip
\noindent
{\it Proof of Lemma~\ref{l_rec_seed_exist}.}
 By Condition~B, for some $\eps >0$ the set
\[
H_0=\{ \omega: \omega (v: |v|\geq 2) \geq \eps\} \cap \supp Q
\]
has positive $Q$-probability.
By the recurrence assumption, the set of $\o$ such that $\omega_0 \in H_0$
and such that
\[
 \Po^y[\mbox{at least one particle hits $0$}] = 1
\]
for any $y\in \A$ has positive $\IP$-probability. We fix $\o'$ in this
set and also in the support of $\IP$.
Then, for any~$\rho<1$ it is possible to choose~$m$ in such a way that
\[
 \min_{y\in \A} {\mathtt P}_{\omega'}^y[
\mbox{at least one particle hits $0$ before time~$m$}] > \rho.
\]
The probability in the above display depends only on the environment inside
the cube $\K_{mL_0}$.
By continuity we can choose neighborhoods
 $H_z\subset\supp Q$ of $\omega_z'$, $z\in\K_{mL_0}$, with  $Q(H_z)>0$,
\[
 \inf_{\o}\min_{y\in \A}
  {\mathtt P}_{\omega}^y[\mbox{at least one particle hits $0$ before time~$m$}] 
                                   > \rho,
\]
where the infimum is taken over all possible environments~$\o$
such that $\omega_z\in H_z$ for all $z\in\K_{mL_0}$.
Due to the boundedness of jumps, for any~$\o\in {\M}^{\Z^d}$ and any $y\in \A$
\begin{eqnarray*}
\lefteqn{\Po^y[\mbox{at least one particle hits $0$ before time~$m$}]} \\
& \leq &
 \PoA{\K_{mL_0}}^y[\mbox{at least one particle hits $0$}].
\end{eqnarray*}
Hence, under $\Po^0$, with probability $\eps$ two particles will be present
at time 1 in $\A$ and otherwise at least one particle, 
each of them having independent
evolution and probability at least $\rho$ to come back to 0 before exiting
$\K_{mL_0}$. By an elementary computation, we see that
\[
\PoA{\K_{mL_0}}^0[ 0 {\rm \ will\ be\ visited\ by\ at\ least\ } i {\rm \
 free\ particles}]
   \geq \sum_{j=i}^2 p_j
\]
with
\begin{equation}
  \label{def15}
p_1=(1-\eps)\rho+2\eps\rho(1-\rho),
p_2=\eps\rho^2, p_3=p_4=\ldots=0.
\end{equation}
It remains only to choose~$m$ large enough to assure that~$\rho$
becomes sufficiently close to~$1$ to guarantee that
the mean $(1-\eps)\rho+2\eps\rho(1-\rho)+2\eps\rho^2$ of
$\p$ defined above is strictly larger than 1.
Then, the $(\K_{mL_0},H)$-seed
constructed in this way is $\p$-recurrent, and it has a positive $\IP$-probability.
\qed

For later purposes, it is useful to emphasize the kind of seed we
constructed above.

\begin{df}
\label{good_seed}
Let $U,W$ be two finite subsets of~$\Z^d$ such that $0\in W$,
$\A+W\subset U$, and let
$\eps,\rho \in (0,1)$.
An $(U,H)$-seed is called $(\eps,\rho,W)$-good, if
\begin{itemize}
\item[(i)] for any $\omega\in H_z$ we have $\omega(v:|v|\geq 2) > \eps$
for all~$z\in W$;
\item[(ii)] for any $\o$ such that $\omega_x\in H_x$, $x\in U$, we have
\[
 \PoA{U}^y\big[\mbox{at least one particle hits~$W$}\big] > \rho
\]
for any $y\in \A+W$;
\item[(iii)] we have $(1-\eps)\rho+2\eps\rho(1-\rho)+2\eps\rho^2 > 1$.
\end{itemize}
In the case~$W=\{0\}$ we say that the seed is $(\eps,\rho)$-good.
\end{df}


At the end of last proof, we just showed that such a seed is
$\p$-recurrent, in the case $W=\{0\}$. It is a simple exercise to extend
the proof  to the case of a general $W$. We state now this useful fact.
\begin{lmm}
\label{l_good_seed}
Any $(U,H)$-seed which is $(\eps,\rho,W)$-good is
also $(\p,W)$-recur\-rent with~$\p$ defined by~(\ref{def15}).
\end{lmm}

Now we finish the proof of Theorem~\ref{trans/rec}. Observe that
if a $\p$-recurrent seed has positive probability
under~$Q$,
then it has positive probability
under~$Q'$
for any~$Q'$ such that $\supp Q = \supp Q'$. (One may invoke the stronger condition
of $Q$ being equivalent to $Q'$, but due to the particular form of the seed
we consider here, the weaker condition of equal support is sufficient.)
An application
of Lemmas~\ref{l_rec_seed_rec} and~\ref{l_rec_seed_exist} concludes the proof.
\qed

\subsection{Proof of Theorem~\ref{suff_rec}}
\label{s_proofrec}
{\bf 1.} We start with the first statement, assuming Condition~E only.
For any $s\in\S^{d-1}$ define
\[
\phi_Q(s) = \sup_{\omega\in\supp Q} \sum_{v\in\V} \omega(v) D(s,v),
\]
with $D$ defined in (\ref{def_D}).
Since $\phi_Q(s)$ is a continuous function of~$s$ and $\S^{d-1}$ is compact,
(\ref{sdf}) implies that
\[
a_0 := \inf_{s\in\S^{d-1}} \phi_Q(s) > 0.
\]
Since $\supp Q$ is closed, for any~$s$ there exists $\omega^{(s)}$
such that
\[
\phi_Q(s) = \sum_{v\in\V} \omega^{(s)}(v) D(s,v).
\]
Moreover, by continuity for any~$s$ we can find $\delta_s>0$ and an open set
$\Gamma_s\subset \M$ with $\omega^{(s)}\in\Gamma_s$ and $Q(\Gamma_s)>0$
such that
\begin{equation}
\label{iop}
\inf_{\substack{s'\in\S^{d-1}:\\ \|s-s'\|<\delta_s}}\inf_{\omega\in\Gamma_s}
             \sum_{v\in\V} \omega(v) D(s',v) > \frac{a_0}{2},
\end{equation}
where $\|\cdot\|$ stands for the Euclidean norm.

Since $\S^{d-1}$ is compact, we can choose $s_1,\ldots,s_m \in\S^{d-1}$ that
generate a finite subcovering of~$\S^{d-1}$ by the sets
$\{s'\in\S^{d-1}:\|s_n-s'\|<\delta_{s_n}\}$, $n=1,\ldots,m$.
For each $n=1,\ldots,m$, it is possible to choose a set
$U_n\subset \{s'\in\S^{d-1}:\|s_n-s'\|<\delta_{s_n}\}$ in such a way that
$U_i\cap U_j = \emptyset$ for $i\neq j$ and $\bigcup_{i=1}^m U_i = \S^{d-1}$.

To prove recurrence, we construct a $(\eps,\rho,W)$-good
$(A,H)$-seed with a supercritical branching inside $W$ and, in
$A\setminus W$, the drift pointing towards $W$ (and so this
seed is a trap).

\begin{figure}
\centering
\includegraphics[width=\textwidth]{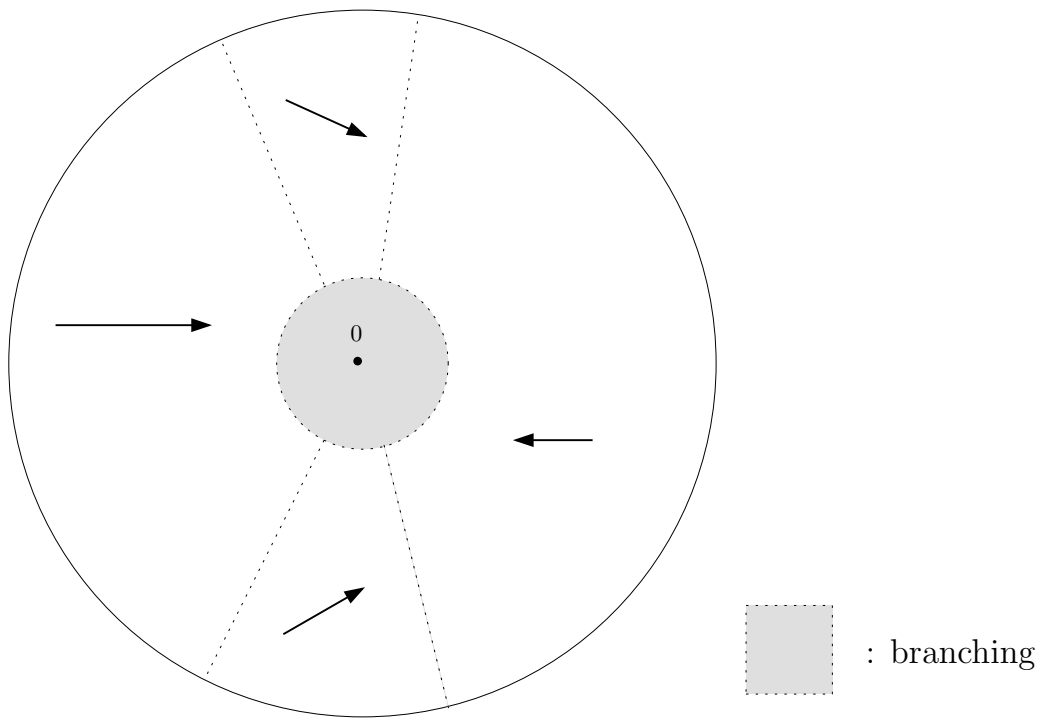}
\caption{Construction of an $(\eps, \rho)$-good 
seed for the branching random walk (defined by~$Q_2$, with $a<1/2$) 
of Example~\ref{exx}}
\label{f_pr-th1-2}
\end{figure}

\begin{itemize}
\item[(i)] similarly to the proof of Theorem~\ref{trans/rec}, we argue that
there exist $\eps>0$ and ${\tilde H}\subset\supp Q$
such that $Q({\tilde H})>0$ and $\omega(v:|v|\geq 2)\geq\eps$
for any $\omega\in {\tilde H}$;
\item[(ii)] take $W=\{y\in\Z^d : \|y\| \leq L_0^2/a_0\}$, and put
$H_z={\tilde H}$
for all $z\in W$;
\item[(iii)] choose $\rho>0$ in such a way that the condition~(iii)
of Definition~\ref{good_seed} holds;
\item[(iv)] choose large enough~$r_2$ in order to guarantee that
       $\rho\leq \frac{r_2-r_1-L_0\sqrt{d}}{r_2-r_1+L_0\sqrt{d}}$, where
       $r_1:=L_0^2/a_0$;
\item[(v)] to complete the definition of the seed, take
       $A=\{y\in\Z^d : \|y\| \leq r_2\}$; it remains to define the
       environment in~$A\setminus W$. It is done in the following way:
       if $z\in A\setminus W$, let~$n_0$ be such that $(-z/\|z\|)\in U_{n_0}$;
       then put $H_z=\Gamma_{s_{n_0}}$, see Figure~\ref{f_pr-th1-2}.
\end{itemize}

To prove that the seed constructed above is indeed $(\eps,\rho,W)$-good,
we construct a random walk~$\xi_n$ that is similar to the example~(iii)
of $r$-extremal induced random walks.
Specifically, suppose that at some moment the random walk~$\xi_n$
is in site~$z\neq 0$ (this does not complicate anything, since we really need
the random walk to be defined only inside $A\setminus W$ 
and $0\notin A\setminus W$).
Generate the offsprings of that particle from~$z$ according to the rules
of the branching random walk; suppose that those offsprings went
to $z+w_1,\ldots,z+w_k$. Let $k_0$ be the number which maximizes the quantity
$(-z/\|z\|)\cdot w_k$; then take $\xi_{n+1} = z+w_{k_0}$. Clearly, for the
random walk constructed in this way,
\begin{equation}
\label{eq_snos}
 \Eo\big((\xi_{n+1}-\xi_n)\cdot (-z/\|z\|) \mid \xi_n=z\big)
             = \sum_{v\in\V} \omega_z(v) D(-z/\|z\|,v).
\end{equation}

Now, we have to bound from below the probability that the random walk~$\xi_n$
starting somewhere from~$W+\A$ hits~$W$ before stepping out from~$A$.
We use ideas which are classical in  Lyapunov functions approach
\cite{FMM}.
To do that, we first recall the following elementary inequality: for
any~$x\geq -1$,
\begin{equation}
\label{koren'}
 \sqrt{1+x} \leq 1 + \frac{x}{2}.
\end{equation}
Let $p_{x,y}$ be the transition probabilities of the random walk~$\xi_n$.
Using~(\ref{koren'}), (\ref{eq_snos}), and~(\ref{iop}), we obtain
\begin{eqnarray}
\Eo(\|\xi_{n+1}\|-\|\xi_n\| \mid \xi_n=z)
& = & \sum_{y\in \A} p_{z,z+y} (\|z+y\|-\|z\|) \nonumber\\
&=& \|z\| \sum_{y\in \A} p_{z,z+y} \Big(\sqrt{1+\frac{2z\cdot y}{\|z\|^2}+
     \frac{\|y\|^2}{\|z\|^2}}-1\Big) \nonumber\\
&\leq & \sum_{y\in \A} p_{z,z+y} \frac{z}{\|z\|}\cdot y + \frac{L_0^2}{2\|z\|}
 \label{1_mom}
 \\
&\leq & -\frac{a_0}{2} + \frac{L_0^2}{2\|z\|} \leq 0 \label{first_mom}
\end{eqnarray}
for all $z\in A\setminus W$. Let~$\tau$ be the first moment when $\xi_n$ leaves
the set $A\setminus W$; the calculation~(\ref{first_mom}) shows that the process
$\|\xi_{n\wedge \tau}\|$ is a (bounded) supermartingale. Denoting by~$\tilde p$
the probability that $\xi_n$ hits~$W$ before stepping out from~$A$ and
starting the walk inside~$W+\A$, we obtain from the Optional Stopping
Theorem that
\[
r_1+L_0\sqrt{d} \geq \Eo\xi_o \geq \Eo\xi_{\tau}
     \geq {\tilde p}(r_1-L_0\sqrt{d}) + (1-{\tilde p})r_2.
\]
So,
\[
 {\tilde p} \geq \frac{r_2-r_1-L_0\sqrt{d}}{r_2-r_1+L_0\sqrt{d}} \geq \rho,
\]
which shows that the $(A,H)$-seed constructed above is
$(\eps,\rho,W)$-good.
With an application of Lemma~\ref{l_good_seed} and of
Lemma~\ref{l_rec_seed_rec}
we conclude the proof of the first
part of Theorem~\ref{suff_rec}.

\medskip
\noindent
{\bf 2.} Now, we  prove that~(\ref{sdf'}) together with Condition~UE
imply recurrence.
The basic idea is the same: we would like to construct a $(A,H)$-seed
that is $(\eps,\rho,W)$-good,
where $W=\{y\in\Z^d : \|y\| \leq r_1\}$, $A=\{y\in\Z^d : \|y\| \leq r_2\}$,
for some $r_1, r_2$ (to be chosen later).
As above, we choose $\eps>0$ in such a way that
there exists ${\tilde H}\subset\supp Q$
such that $Q({\tilde H})>0$ and $\omega(v:|v|\geq 2)\geq\eps$
for any $\omega\in {\tilde H}$.
Then, choose $\rho>0$ in such a way that the condition~(iii)
of Definition~\ref{good_seed} holds. To define the seed inside~$W$,
put $H_z={\tilde H}$ for all $z\in W$. Now, it remains to specify $r_1, r_2$
and to define the seed in $A\setminus W$. To do that, we need some
preparations. For any {\it possible\/} environment inside $A\setminus W$
(i.e., for any~$\o$ such that $\omega_z\in\supp Q$ for all $z\in A\setminus W$)
we keep the same definition of the random walk~$\xi_n$; by
Condition~UE, there exists $\gamma>0$ such that for all~$x\neq 0$
\begin{equation}
\label{eq_C'}
 \Eo\big((\|\xi_{n+1}\|-\|\xi_n\|)^2\1{\|\xi_{n+1}\|\leq\|\xi_n\|} 
            \mid \xi_n=x\big)
      \geq \gamma.
\end{equation}
With that~$\gamma$, we successively choose $\alpha>0$ such that
\begin{equation}
\label{eq_choose_alpha}
 \gamma(\alpha+1) > L_0^2,
\end{equation}
then $r_2>r_1>L_0$ with
\begin{equation}
\label{eq_choose_r12}
 \frac{(r_1+L_0)^{-\alpha}-r_2^{-\alpha}}{(r_1-L_0)^{-\alpha}-r_2^{-\alpha}}
      > \rho,
\end{equation}
and finally,  $\eps'>0$ such that
\begin{equation}
\label{eq_choose_eps''}
\eps' < \frac{\gamma(\alpha+1)-L_0^2}{2r_2}.
\end{equation}
Now we define the seed on the set $A\setminus W$ in the following way.
When~(\ref{sdf'}) holds, analogously to the first part of the proof of this
theorem, for any~$s\in\S^{d-1}$ we can find $\delta'_s>0$ and an open set
$\Gamma'_s\subset \M$ with $\omega^{(s)}\in\Gamma'_s$ and $Q(\Gamma'_s)>0$
such that
\begin{equation}
\label{iop'}
\inf_{\substack{s'\in\S^{d-1}:\\ \|s-s'\|<\delta'_s}}\inf_{\omega\in\Gamma'_s}
             \sum_{v\in\V} \omega(v) D(s',v) > - \eps'.
\end{equation}
Similarly to what has been done before,
we choose $s'_1,\ldots,s'_m \in\S^{d-1}$ that
generate a finite subcovering of~$\S^{d-1}$ by the sets
$\{s\in\S^{d-1}:\|s'_n-s\|<\delta_{s'_n}\}$, $n=1,\ldots,m$.
For each $n=1,\ldots,m$, we choose a set
$U'_n\subset \{s\in\S^{d-1}:\|s'_n-s\|<\delta_{s'_n}\}$ in such a way that
$U'_i\cap U'_j = \emptyset$ for $i\neq j$ and $\bigcup_{i=1}^m U'_i = \S^{d-1}$.
Now, if $z\in A\setminus W$, let~$n_1$ be such that $(-z/\|z\|)\in U'_{n_1}$;
       then put $H_z=\Gamma_{s'_{n_1}}$.

To prove that the $(A,H)$-seed is indeed $(\eps,\rho,W)$-good,
we have to verify the condition~(ii) of Definition~\ref{good_seed}.
First, it elementary to see that the following inequality holds:
for any $x\geq -1$, $\alpha>0$
\begin{equation}
\label{eq_***}
 (1+x)^{-\alpha} \geq 1-\alpha x +
            \frac{\alpha(\alpha+1)}{2}x^2\1{x\leq 0}.
\end{equation}
Keeping the notation~$\tau$ from the proof of the first part of the theorem,
we are aiming to prove that $\|\xi_{n\wedge\tau}\|^{-\alpha}$ is
a submartingale. Indeed, using (\ref{eq_***}), (\ref{1_mom}),
(\ref{eq_C'}) and (\ref{eq_choose_eps''}), we obtain for $z\in A\setminus W$
\begin{eqnarray*}
\lefteqn{\Eo(\|\xi_{n+1}\|^{-\alpha}-\|\xi_n\|^{-\alpha} \mid \xi_n=z)} \\
 &=& \|z\|^{-\alpha}
 \Eo\Big(\Big(1+\frac{\|\xi_{n+1}\|-\|\xi_n\|}{\|\xi_n\|}\Big)^{-\alpha} - 1
    \,\Big|\, \xi_n=z\Big) \\
 & \geq & -\alpha \|z\|^{-\alpha-1}\Eo(\|\xi_{n+1}\|-\|\xi_n\| \mid \xi_n=z)\\
 && {}+ \frac{\alpha(\alpha+1)}{2}\|z\|^{-\alpha-2}
      \Eo((\|\xi_{n+1}\|-\|\xi_n\|)^2\1{\|\xi_{n+1}\|\leq\|\xi_n\|} \mid \xi_n=z)\\
 &\geq & \alpha \|z\|^{-\alpha-1}\Big(-\eps'-\frac{L_0^2}{2\|z\|}+
                   \frac{(\alpha+1)\gamma}{2\|z\|}\Big) \\
 & > & 0.
\end{eqnarray*}
Then, by using the Optional Stopping Theorem again, we obtain that
the probability that $\xi_n$ hits~$W$ before stepping out from~$A$ and supposing
that its starting point belongs to~$W+\A$ is at least
\[
 \frac{(r_1+L_0)^{-\alpha}-r_2^{-\alpha}}{(r_1-L_0)^{-\alpha}-r_2^{-\alpha}},
\]
so, recalling~(\ref{eq_choose_r12}),
we see that the condition~(ii) of Definition~\ref{good_seed} holds.
We finish the proof of the second part of Theorem~\ref{suff_rec}
by applying Lemma~\ref{l_good_seed}.    \qed

\subsection{Proof of Theorem~\ref{suff_tr}}

Due to Condition~B, there exists $\omega\in\supp Q$
such that $\sum_{y\in \A}\mu_y^\omega > 1$. Hence, Condition~L
cannot be satisfied with $\lambda=1$.
Moreover, if Condition~L holds
for~$\lambda$ and~$s$, it holds also for~$\lambda^{-1}$ and~$(-s)$.
So, we can suppose that $\lambda\in (0,1)$ without loss of
generality.

Denote the half-space $Y_s = \{y\in\Z^d: y\cdot s\leq 0\}$.
Take an arbitrary starting point $z\notin Y_s$ and define
\[
 F_n^z = \sum_{y\in\Z^d}\eta_n^z(y) \lambda^{y\cdot s}.
\]
{\bf 1.} Let us also modify the environment in such a way that any particle which
enters~$Y_s$ neither moves nor branches anymore.
By~(\ref{eq_C_L}) it is straightforward to obtain that
the process $(F_n^z, n=0,1,2,\ldots)$ is a supermartingale:
\begin{eqnarray*}
 \Eo^z\big(F_{n+1}^z | \eta_1^z(\cdot), \ldots, \eta_n^z(\cdot)\big)
&=&
 \sum_{x \in Y_s} \eta_n^z(x) \lambda^{x\cdot s}+
 \sum_{x \notin Y_s} \eta_n^z(x) \lambda^{x\cdot s} \times
\sum_{y \in \A}
\mu_y^{\omega_x} \lambda^{y\cdot s}\\
&\leq & F_n^z.
\end{eqnarray*}
Since it is also nonnegative, it converges a.s.\ as $n\to\infty$ to
some random variable~$F_\infty$. By Fatou's lemma,
\[
 \Eo^z F_\infty \leq \Eo^z F_0 = \lambda^{z\cdot s} < 1 ,
\]
for $z \notin Y_s$.
On the other hand, any particle stuck in~$Y_s$
contributes at least one unit to~$F$. That shows that with positive probability
the branching random walk will not enter to~$Y_s$, so
the proof of
the first part of
Theorem~\ref{suff_tr} is finished.
\medskip

\noindent
{\bf 2.} We no longer make $Y_s$ absorbing. Note that $F_n^z$ is still a
supermartingale, and has an a.s.\ limit (for all $\o$).
Let $k \geq 1$.
First of all, let us show how to prove the result when Condition~UE holds.
Under Condition~UE, each time a particle enters the half-space $Y_s$
it has $\Po$-probability larger than $\eps_0^{dk}$ to enter
  $Y_s^{(k)}= \{y\in\Z^d: y\cdot s\leq - k\}$.
By the strong Markov property, an infinite number of particles will hit
$Y_s^{(k)}$
a.s.\ on the set where the number of visit of the \BRW \ to~$Y_s$ is
infinite. We will then have, on this set, $\limsup_n F_n^z \geq
\lambda^{-k}$ for all $k$ (recall that $\lambda<1$). 
Since $F_n^z$ has a finite limit, this
shows that the number of visits to $Y_s$ is finite.

Now, we explain what to do when only Condition~E holds.
By the previous argument, it would be enough to prove that,
on the event that~$Y_s$ is visited infinitely often,
for any~$k$, the set $Y_s^{(k)}$ is visited infinitely many times $\Po^z$-a.s.
Define $H_s^{(k)} = Y_s\setminus Y_s^{(k)}$.
Suppose, without restriction of generality, that $\|s\|_1=1$.
For any $z\in H_s^{(k)}$ define
\[
 g_z^{(k)}(\o) = \Po^z[\mbox{at time $k$ there is at 
                      least one particle in $Y_s^{(k)}$}].
\]
It is elementary to obtain that
\begin{itemize}
\item[(i)] $g_u^{(k)}(\o)$ and $g_v^{(k)}(\o)$ are 
independent if $\|u-v\|_\infty > 2kL_0$, and
\item[(ii)] there exists $h_k>0$ such that $\IP[g_z^{(k)}(\o)>h_k]\geq 1/2$, 
uniformly in $z\in H_s^{(k)}$.
\end{itemize}

Now, suppose that $Y_s$ was visited an infinite number of times, and suppose also
that the number of visits to $H_s^{(k)}$ is also infinite (because otherwise,
automatically, $Y_s^{(k)}$ is visited infinitely many times). 
Let $z_1,z_2,z_3,\ldots$
be the locations of those visits. Using (i) and~(ii),
we can extract
an infinite subsequence $i_1<i_2<i_3<\ldots$ such that
$g_{z_{i_j}}^{(k)}(\o) > h_k$, for all~$j$. Similarly to the previous argument,
we obtain that in this case $Y_s^{(k)}$ will be visited infinitely often,
which leads to a contradiction with the existence of a finite limit
for $F_n^z$.
\qed

\section{Proof of Theorems~\ref{t_upbound_ht} and~\ref{t_lobound_ht}}
\label{s_proofs_ht}

\noindent
{\it Proof of Theorem~\ref{t_upbound_ht}.}
 Roughly, the idea is as follows: by recurrence we know
that there are $\p$-recurrent (in fact, even $(\eps,\rho)$-good)
seeds, each of them supporting a
supercritical Galton-Watson process. To prove~(\ref{eq_hit_time}) it suffices
essentially to control the time to reach a large
enough quantity of these seeds.

By Lemma~\ref{l_rec_seed_exist}, there exist $n_0,\eps,\rho>0$
and a collection $H=(H_z\subset\M, z\in\K_{n_0})$
having positive $\IP$-probability,
such that the $(\K_{n_0},H)$-seed is $(\eps,\rho)$-good
(in the proof of Lemma~\ref{l_rec_seed_exist}, we indeed constructed
such a seed).
Moreover, it is straightforward to see that there exists~$t_0$ such that
\begin{eqnarray}
\IP\Big[ \mbox{the $(\K_{n_0},H)$-seed is $(\eps,\rho)$-good and for any
         $y\in\A$}~~~~~~~~~~\nonumber\\
\PoA{\K_{n_0}}^y[\mbox{at least one particle hits~$0$ before
time~$t_0$}]>\rho\Big]
  & > & 0.\phantom{***} \label{improve_good}
\end{eqnarray}
For any~$\o$, define the random subset $S_{\o}$
of the lattice with spacing $2n_0+1$
\begin{eqnarray*}
S_{\o} & = &\{z\in (2n_0+1)\Z^d : z \mbox{ is the center of $(\K_{n_0},H)$-seed} \\
&&~~~~~~~~~~~~~~~~~~~~~~\mbox{which is $(\eps,\rho)$-good
              and satisfies~(\ref{improve_good})}\}.
\end{eqnarray*}

We need to consider two cases separately, $d\geq 2$ and $d=1$.

\medskip
\noindent
{\bf Case $d\geq 2$.} 
Consider the event
\begin{eqnarray*}
 M_n &=& \{\forall y\in\K_{L_0 n\ln^{-1} n}
 \exists z\in S_{\o} :  \|y-z\|_{\infty} \leq \alpha\ln n\},
\end{eqnarray*}
the (small enough) constant~$\alpha$ will be chosen later.
We will use the bound
\begin{equation}
\label{ann_prob}
 \PA[T(0,x_0)>n] \leq \sup_{\o\in M_n} \Po^0[T(0,x_0)>n]
                              + \IP[M_n^c].
\end{equation}
Let us begin by estimating the second term 
in the right-hand side of~(\ref{ann_prob}).
We have
\begin{eqnarray}
\IP[M_n^c] &=& \IP[\exists y\in\K_{L_0 n\ln^{-1} n} :
       (y+\K_{\alpha\ln n})\cap S_{\o}=\emptyset] \nonumber\\
       &\leq & |\K_{L_0 n\ln^{-1} n}| \IP[\K_{\alpha\ln n}\cap S_{\o}=\emptyset].
       \label{ocenka_compl}
\end{eqnarray}
The point is that
the events $\{x\in S_{\o}\}$ and $\{y\in S_{\o}\}$ are independent
for any $x,y\in (2n_0+1)\Z^d$, $x\neq y$.  Denoting the left-hand side
of~(\ref{improve_good}) by
$p_0=\IP[0\in  S_{\o}]>0$, we obtain
\[
\IP[\K_{\alpha\ln n}\cap S_{\o}=\emptyset] \leq
                    (1-p_0)^{\frac{\alpha^d\ln^d n}{(2n_0+1)^d}},
\]
so, from~(\ref{ocenka_compl}),
\begin{eqnarray}
\IP[M_n^c] & \leq & L_0^d n^d \ln^{-d}n
   \exp\Big\{-\frac{\alpha^d\ln(1-p_0)^{-1}}{(2n_0+1)^d}\ln^d n\Big\} \nonumber\\
 & \leq & \exp\{-C_1\ln^d n\} \label{term_2}
\end{eqnarray}
for some $C_1>0$ and for all~$n$ large enough.

Now, we estimate the first term in the right-hand side of~(\ref{ann_prob}).
Let~$\xi_n$ be the uniform induced random walk in random environment,
cf.\ example~(i) in Section~\ref{sec:induced}.
By Condition~UE, this random walk will be uniformly elliptic as well,
in the sense that for any~$x\in\Z^d$ and any $\omega\in\supp Q$
\begin{equation}
\label{ell_induced}
 \Po^x[\xi_1 = x+e] \geq \eps_1 > 0
\end{equation}
for all $e\in \{\pm e_i, i=1,\ldots,d\}$
with a new constant $\eps_1=\eps_0 (2L_0+1)^{-1}$.
By~(\ref{ell_induced}),
we have that for arbitrary $\o\in M_n$ and any~$m$
with $0 \leq m \leq n\ln^{-1}n  - d\alpha\ln n$,
\begin{equation}
\label{enter_S}
 \Po^0\big[\{\xi_m,\ldots,\xi_{m+d\alpha\ln n}\} \cap S_{\o} \neq
   \emptyset
\;|\; \xi_0,\ldots,\xi_{m-1}\big]
               \geq \eps_1^{d\alpha\ln n}.
\end{equation}

Define
\begin{equation}
\label{def_tau}
 \tau = \inf\Big\{m:\sum_{i=0}^m\1{\xi_i\in S_{\o}} \geq \ln^d n\Big\},
\end{equation}
i.e., $\tau$ is the moment when the random walk~$\xi$ hits the
set~$S_{\o}$ for the $\lceil\ln^d n\rceil$-th time.
Also, let us recall the Chernoff's bound for the binomial distribution:
if $S_k$ is a binomial ${\cal B}(n,p)$ random variable,
for any~$k$ and $a$   with $0<a<p<1$, we have
\begin{equation}
\label{Chernoff_a<p}
\PA\Big[\frac{S_k}{k} \leq a\Big] \leq \exp\{-k U(a,p)\},
\end{equation}
where
\[
 U(a,p) = a\ln\frac{a}{p} + (1-a)\ln\frac{1-a}{1-p} > 0.
\]
Now, divide the time interval $[0,n\ln^{-1}n]$ into
$(d\alpha)^{-1}n\ln^{-d}n$ subintervals of length $d\alpha\ln n$.
Fix the constant~$\alpha$ in such a way that $d\alpha \ln \eps_1^{-1}<1/2$.
Use Markov property for $\xi$ under $\Po^0$,
the inequality~(\ref{enter_S}), and~(\ref{Chernoff_a<p})
with $p=\eps_1^{d\alpha\ln n}$, $k=(d\alpha)^{-1}n\ln^{-d}n$,
$a=d\alpha n^{-1}\ln^{2d}n$ (and an elementary computation shows  
that $U(a,p)$ is of order $n^{-d\alpha \ln \eps_1^{-1}}$) 
to obtain that for some $C_2,C_3>0$
\begin{eqnarray}
\Po^0[\tau \leq n/3] & \geq & \Po^0[\tau \leq n\ln^{-1}n] \nonumber\\
&\geq & \Po^0\Big[\sum_{i=1}^k {\bf 1}_{\{\xi_j; (i-1)d\alpha\ln n <j\leq
i\alpha d\ln n\} \bigcap {S_{\o}} {\neq \emptyset} } \geq ka \Big]
\nonumber\\
 & \geq & 1-\exp\big\{-C_2(d\alpha)^{-1}n^{1-d\alpha \ln \eps_1^{-1}} 
                 \ln^{-d} n   \big\}\nonumber\\
  & \geq & 1-\exp\{-C_3 n^{1/2}\} \label{ocenka_tau}
\end{eqnarray}
for any $\o\in M_n$ (supposing that~$n$ is large enough so that
$a<p$).

Now, we show that each time the random walk~$\xi$ passes through the points
of~$S_{\o}$ it gives rise to a supercritical Galton-Watson process,
and that on the set $\{\tau \leq n/3\}$, about $\ln^d n$ such independent
Galton-Watson processes will be started before time $n/3$.
Indeed, analogously to the proof of Lemma~\ref{l_good_seed},
if we have a particle in the center of the seed, its direct offsprings
in this Galton-Watson process are those descendants
(in the branching random walk restricted on the seed) that pass through
the center not later than~$t_0$. (Actually, we must take
this Galton-Watson process
independent of the random walk~$\xi$, so when~$\xi$ passes through the seed, we
cannot use the corresponding particle in the Galton-Watson process.
This, however, does not spoil anything, because with uniformly
positive probability
another particle will be generated somewhere in the set $x+\A$
--with~$x$ the center--
at that moment, so
it can be used to start the Galton-Watson process.) By construction,
this Galton-Watson process is ``uniformly'' supercritical, so
there exists $p_1>0$ such that with
probability at least~$p_1$ in the $[n/3t_0]$-th generation
of the process the number of particles will be at least~$C_4\alpha_1^n$,
for some $C_4>0$, $\alpha_1>1$. So, since the real time between the
generations is at most~$t_0$, this means that for any $x\in S_{\o}$
\begin{eqnarray}
\Po^x[\mbox{the seed in $x$ generates at least~$C_4\alpha_1^n$ free
particles}\phantom{*} \nonumber\\
\mbox{ before time $n/3$}] &>& p_1.\phantom{***}
 \label{eq_blow}
\end{eqnarray}

Now, by~(\ref{eq_blow}) we have
\begin{multline}
\Po^0\Big[\mbox{at least one seed generates at least~$C_4\alpha_1^n$
free particles}\\
\mbox{ before time $(2n/3)$} \mid \tau<n/3\Big]
\geq
1-(1-p_1)^{\ln^d n}.
\label{eq_atleastone}
\end{multline}
Consider those $C_4\alpha_1^n$ free particles. By Condition~UE,
any
descendants of each one will hit~$x_0$ by the time $2L_0dn\ln^{-1}n<n/3$
with probability at least $\eps_0^{2L_0dn\ln^{-1}n}$, so at least one
particle will hit~$x_0$ with probability at least
\[
1-(1-\eps_0^{2L_0dn\ln^{-1}n})^{C_4\alpha_1^n} \geq
 1-\exp\Big\{-C_4\exp\{n[\ln\alpha_1-2L_0d\ln\eps_0^{-1}\ln^{-1}n]\}\Big\}
\]
where the
square parentheses is positive for large enough~$n$.
Taking into account~(\ref{ocenka_tau}) and~(\ref{eq_atleastone}),
we then obtain that
for any $\o\in M_n$
\begin{equation}
\label{quench_prob}
\Po^0[T(0,x_0) > n] \leq e^{-C_5\ln^d n}
\end{equation}
for some $C_5>0$ and all~$n$ large enough.
We plug now~(\ref{term_2}) and~(\ref{quench_prob}) into~(\ref{ann_prob})
to conclude the proof of Theorem~\ref{t_upbound_ht} in the case $d\geq 2$.

\medskip
\noindent
{\bf Case $d=1$.}
Now, we prove Theorem~\ref{t_upbound_ht} in dimension~$1$. For~$d=1$
the above approach fails, because if~$\alpha$ is small, then~(\ref{term_2})
will not work, if~$\alpha$ is large, then we would have 
problems with~(\ref{ocenka_tau}), and it is not always possible
to find a value of~$\alpha$ such that both inequalities would work.

First, we do the proof assuming that $L_0=1$, i.e., $\A$ is either
$\{-1,1\}$ or $\{-1,0,1\}$. Analogously to the proof for higher
dimensions, if we prove that, on the set of environments of
$\IP$-probability at least $1-e^{-C_1\ln n}$, the initial particle
hits a quantity of logarithmic order of good seeds from~$S_\o$,
we are done. To this end, note that, on the time interval of length
$\frac{\ln n}{2\ln\eps_0^{-1}}$ a single particle (even if it does
not generate any offsprings) covers a space interval of the same
length with probability at least 
\[
  \eps_0^{\frac{\ln n}{2\ln\eps_0^{-1}}} = n^{-1/2}.
\]
So, by time $n\ln^{-1}n$, with large (of at least stretched exponential order) 
probability there is an interval of 
length $\frac{\ln n}{2\ln\eps_0^{-1}}$ containing~$0$ such that
all sites from there are visited.

Analogously to the proof for higher dimensions, consider the set
\begin{eqnarray*}
 M_n^{(1)} &=& \{\text{the number of good seeds 
                    from~$S_\o$ in all the intervals}\\
&& ~~~~~~~~~~~\text{of length 
$\textstyle\frac{\ln n}{2\ln\eps_0^{-1}}$ containing~$0$ is at least
$C_2\ln n$}\},
\end{eqnarray*}
which corresponds to the set of ``good'' environments. Since~$S_\o$
has a positive density, we can choose small enough~$C_2$ in such a
way that
\[
  \IP[(M_n^{(1)})^c] \leq e^{-C_3\ln n}
\]
for some~$C_3$. Now, on~$M_n^{(1)}$, with probability at least
$1-e^{-C_4n^{C_5}}$ by time $n\ln^{-1}n$ at least $C_2\ln n$ good seeds 
from~$S_\o$ will be visited. The rest of the proof is completely analogous
to the proof for $d\geq 2$.

Let us explain how to proceed in the case of a general~$L_0$.
In the above argument the fact $L_0=1$ was used only for the
following purpose: if we know that a particle crossed a (space)
interval, then we are sure that all the good seeds that might
be there were visited. For a general~$L_0$, instead of $(\eps,\rho)$-good
seeds of~$S_\o$, use $(\eps,\rho,W)$-good seeds with $W=\{0,1,\ldots,L_0-1\}$,
so that particles cannot jump over the translates of this~$W$.
Indeed, it is clear that a recurrent
branching random walk generates $(\eps,\rho,W)$-good seeds for any
finite~$W$). So, proof of Theorem~\ref{t_upbound_ht} is concluded.
\qed

\medskip
\noindent
{\it Proof of Theorem~\ref{t_lobound_ht}.}
The method of the proof is very similar to the construction of 
Example~\ref{ex_infexp}. Roughly speaking, for given~$n$ and~$x$,
we create a (rather improbable)
environment that has a trap near the origin; for such an environment
with a good probability the event $\{T(0,x)>n\}$ occurs.

Now, let us work out the details.
Rather than doing the proof for $T(0,x)$ with a general $x\in\Z^d$, 
we use $x=e_1$, the general case being completely analogous.
Suppose that the origin belongs to the 
interior of the convex hull of $\{\Delta_\omega :
\omega\in \G\cap \supp Q\}$. Then, analogously to the proof
of Theorem~\ref{suff_rec} (see Section~\ref{s_proofrec}),
one can split the sphere $\S^{d-1}$ into a finite number
(say, $m_0$) of non-intersecting subsets ${\hat U}_1,\ldots, {\hat U}_{m_0}$
and find a finite collection 
${\hat \Gamma}_1,\ldots, {\hat \Gamma}_{m_0}\subset \G$
having the following properties:
\begin{itemize}
\item[(i)] there exists~$p_1>0$ such that $Q({\hat \Gamma}_i)>p_1$,
\item[(ii)] there exists~$a_1>0$ such that
for any $z\in {\hat U}_i$ and any $\omega\in {\hat \Gamma}_i$
we have $z\cdot \Delta_\omega < -a_1$,
\end{itemize}
for all $i=1,\ldots,m_0$.

Take $A=\{y\in\Z^d:\|y\|\leq u\ln n\}$, 
where $u$ is a (large) constant to be chosen later. Consider the
$(A,H)$-seed with $H_x,x\in A$ defined as follows.
First, put $H_0=\G$; for $x\neq 0$, let $i_0$ be such that
$\frac{x}{\|x\|}\in {\hat U}_{i_0}$ (note that~$i_0$ is uniquely
defined), then put $H_x={\hat \Gamma}_{i_0}$.
Clearly
\begin{equation}
\label{prob_trap}
 \IP[\text{there is $(A,H)$-seed in $y$}] \geq p_1^{(2u)^d \ln^d n}.
\end{equation}

Note that for any possible environment inside the $(A,H)$-seed there
is no branching. This means that the process restricted on~$A$
is a random walk (without branching), which will be denoted by~$\xi_n$.

Analogously to~(\ref{first_mom}), we can prove that there exist
$a_2>0,C_1\geq 0$ such that 
\begin{equation}
\label{snos_k_0}
\Eo(\|\xi_{n+1}\|-\|\xi_n\| \mid \xi_n=z) < -a_2
\end{equation}
for all $z\in A\setminus \{y:\|y\|\leq C_1\}$, provided there is
an $(A,H)$-seed in~$0$. Let $\tilde\tau$ be the hitting time
of the set $(\Z^d\setminus A)\cup \{y:\|y\|\leq C_1\}$ by~$\xi_n$.
Next, we prove that, when $a_3>0$ is small enough, the 
process $e^{a_3\|\xi_{n\wedge{\tilde\tau}}\|}$ is a supermartingale.
Indeed, first, note that there exist $C_2,C_3>0$ such that
\begin{equation}
\label{elem_ner}
   e^x < 1+x+C_2x^2
\end{equation}
when $|x|<C_3$. We can choose~$a_3$ small enough so that
$\big|\|\xi_{n+1}\|-\|\xi_n\|\big|<C_3/a_3$ 
a.s. From~(\ref{elem_ner}) we obtain
\begin{eqnarray*}
\Eo(e^{a_3\|\xi_{n+1}\|}-e^{a_3\|\xi_n\|}\mid \xi_n=z)
 & =& e^{a_3\|z\|} \Eo(e^{a_3(\|\xi_{n+1}\|-\|\xi_n\|)}-1\mid \xi_n=z)\\
 &\leq & e^{a_3\|z\|} (-a_3a_2 + C_2a_3^2 L_0^2)\\
&<& 0
\end{eqnarray*}
if $a_3$ is small enough, so $e^{a_3\|\xi_{n\wedge{\tilde\tau}}\|}$ is 
indeed a supermartingale.

Now, we need to make two observations concerning the exit probabilities. 
First, consider any~$y$ such that $C_1+1 \leq \|y\| < C_1+2$.
If~${\hat p}_1$ is the probability that, starting from~$y$,
the random walk~$\xi_n$ hits the set $\Z^d\setminus A$
before the set $\{y:\|y\|\leq C_1\}$, then it is straightforward
to obtain from the Optional Stopping Theorem that
\begin{equation}
\label{oc_naruzhu}
 {\hat p}_1 \leq \frac{e^{a_3(C_1+2)}}{n^{a_3u}}.
\end{equation}
\begin{figure}
\centering
\includegraphics[width=0.71\textwidth]{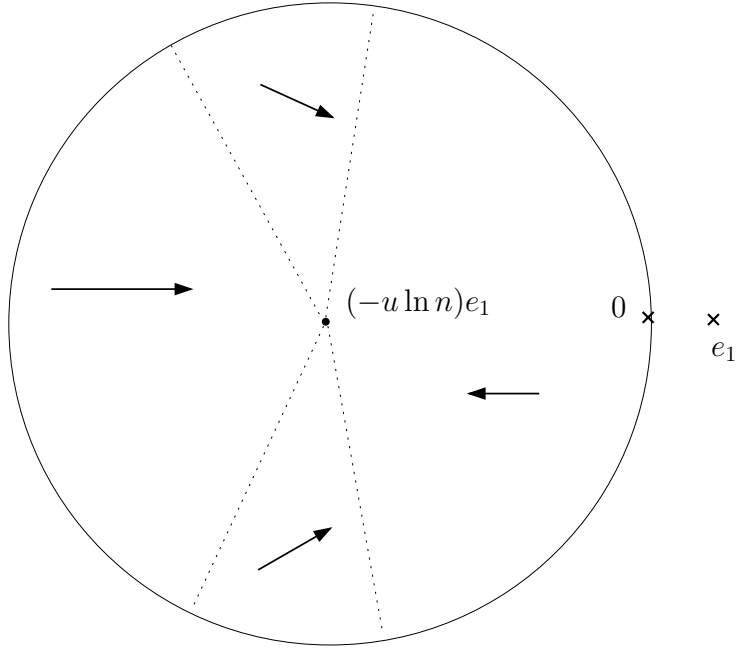}
\caption{Construction of a trap}
\label{f_piege}
\end{figure}
Secondly, suppose now that the random walk~$\xi_n$ starts
from a point~$y$ with $\|y\|=u\ln n$ (i.e., on the boundary of~$A$).
Analogously, using the Optional Stopping Theorem, one can show
that, with probability bounded away from~$0$, the random walk
hits the set $\{y:\|y\|\leq C_1\}$ before stepping out of~$A$.
Now, suppose that $u>a_3^{-1}$, and that there is an $(A,H)$-seed centered
at $(-u\ln n)e_1$ (i.e., touching the origin, cf.\ Figure~\ref{f_piege}). 
Using the previous observation together with~(\ref{oc_naruzhu}),
one can obtain that, with probability bounded away from~$0$,
the particle will go to the set $(-u\ln n)e_1+A$ and will stay there
until time~$n$ (without generating any other particles, since
there is no branching in the sites of $(A,H)$-seed).
So, by~(\ref{prob_trap}),
\[
\PA[T(0,e_1)>n] \geq e^{-C_4\ln^d n},
\]
thus completing the proof of Theorem~\ref{t_lobound_ht}. \qed

\section{Proof of Theorem~\ref{shape_t}}
\label{pf-shape_th}
We prove this theorem separately for two cases: $d\geq 2$
and $d=1$. There are, essentially, two reasons for splitting
the proof into these two cases. First, as usual, in
dimension~$1$ we have to care about only one (well, in fact, two)
directions of growth, while for $d\geq 2$ there are infinitely
many possible directions. So, one may think that the
proof for $d=1$ should be easy when compared to the
proof for $d\geq 2$. For the majority of growth models this is
indeed true, but not for the model of the present paper.
This comes from Theorems~\ref{t_upbound_ht}, \ref{t_lobound_ht},
and Example~\ref{ex_infexp}: 
recurrence implies that the annealed expectation of the hitting time 
is finite only for $d\geq 2$, but not for $d=1$.

\subsection{Case $d\geq 2$}
\label{s_d_geq2}
First, we need to show that the sets of interest grow
at least linearly.
\begin{lmm}
\label{hit_time}
Suppose that $d\geq 2$ and the branching random walk in random environment
is recurrent and Condition~UE holds.
Then
\begin{itemize}
\item[(i)] There exist $\delta_0, \theta_0>0$ such that
\begin{equation}
\label{eq_small_ball'}
\PA[\K_{\delta_0 n}\subset B^0_n] \geq 1 - \exp\{-\theta_0\ln^d n\}
\end{equation}
for all
$n$ sufficiently large.
\item[(ii)] Suppose, in addition, that Condition~A holds.
Then  there exist $\delta_1, \theta_1>0$ such that
\begin{equation}
\label{eq_small_ball}
\PA[\K_{\delta_1 n}\subset {\tilde B}^0_n] \geq 1 -
                                \exp\{- \theta_1\ln^d n\}
\end{equation}
for all
$n$ sufficiently large.
\end{itemize}
\end{lmm}

\noindent
{\it Proof of Lemma \ref{hit_time}.}
We will use the notations from the proof of Theorem~\ref{t_upbound_ht}.

\medskip
\noindent
{\bf Step 1.}
Let us prove part~(ii) first. To do that, we need to examine in more
detail the supercritical Galton-Watson process arising in seeds centered
in the points of~$S_{\o}$. Specifically, we need more information about
how (conditioned on survival) the particles of that process are
distributed in time. As we have seen before, in that Galton-Watson process
a particle has~$1$ offspring with
probability $(1-\eps)\rho+2\eps\rho(1-\rho)$, $2$ offsprings
with probability $\eps\rho^2$ and~$0$ offsprings
with the remaining probability,
and the parameters $\eps,\rho$ are such
that $(1-\eps)\rho+2\eps\rho(1-\rho) + 2\eps\rho^2>1$,
so the process is uniformly supercritical.
Moreover, the real time interval between
the particle and its offspring(s) is not
larger than~$t_0$; the exact distribution
of that time interval, is, however, unknown.
So, let us suppose that if a particle
has~$1$ offspring, than it reappears in the center of the  seed after~$k$
time units with probability~$q_k$, $k=1,\ldots,t_0$, and if a particle has~$2$
offsprings, then they reenter the center after~$i$ and~$j$ time units
with probability~$q_{i,j}$, $i,j=1,\ldots,t_0$, $i\leq j$
(we have $\sum_k q_k = 1$, and $\sum_{i,j}q_{i,j} = 1$).

Supposing for a moment that the Galton-Watson process starts from one particle
at time~$0$, denote by $\zeta(n)$ the number of particles
of that process at time~$n$. We are going to prove that, conditioned
on survival, $\sum_{i=nt_0}^{(n+1)t_0}\zeta(i)$ grows rapidly in~$n$.
For that, we construct two processes $Z_n^i, {\hat Z}_n^i$,
$n=0,1,2,\ldots$, $i=0,\ldots,t_0-1$.
We start by defining ${\hat Z}_0^i=0$ for
all~$i$, $Z_0^i=0$ for $i=1,\ldots,t_0-1$, and $Z_0^0=1$. Inductively,
suppose that the processes $Z,{\hat Z}$ are constructed up to~$n$.
Suppose for example that $Z_n^{i_0}=a>0$; this means that there are~$a$
particles of~$Z$ in the center of the seed at time~$i_0+nt_0$.
For each of those~$a$ particles, do the following:
\begin{itemize}
\item let it generate its offsprings according to the rules of
the Galton-Watson process; those offsprings reenter the center either
during the time interval $[nt_0,(n+1)t_0)$, or during $[(n+1)t_0,(n+2)t_0)$;
\item for those offsprings that appeared in the center of the seed
during the interval $[nt_0,(n+1)t_0)$, repeat the above step.
\end{itemize}
Doing that, we obtain a cloud of free particles
(again, in the sense that one cannot be descendant of another) in
the interval $[(n+1)t_0,(n+2)t_0)$. Fix a parameter~$h>0$
and declare each of those particles to be of type~$1$ with probability~$1-h$
and to be of type~$2$ with probability~$h$, independently.
Repeat the same procedure for all $i_0\in\{0,\ldots,t_0-1\}$
(note that the particles from~${\hat Z}$ are not used in this
construction).
Then, define $Z_{n+1}^i$ to be the number of type~$1$ particles
at the moment $i+(n+1)t_0$, and
${\hat Z}_{n+1}^i$ to be the number of type~$2$ particles
at the same moment, $i=0,\ldots,t_0-1$. Then, $Z_n=(Z_n^0,\ldots,Z_n^{t_0-1})$
is a multitype branching process with $t_0$ types.
Furthermore, it is straightforward to
see that if~$h$ is small enough, then the mean number of particles
(of all types) generated by a particle from~$Z_n$ is greater than~$1$, so
that process is supercritical (this follows from, e.g., Theorem~2 of Section~3
of Chapter~V of~\cite{AN}, noting also that if, for a nonnegative matrix,
the sum of entries is strictly greater than~$1$ for each row, then the maximum
eigenvalue of that matrix is strictly greater than~$1$).
That is, with positive probability
the size of $n$-th generation of~$Z$ grows exponentially in~$n$.
  From that it is quite elementary to obtain that there exist
constants $\gamma_2,p_2>0$, $\alpha_2>1$ (depending on $q_i$-s and
$q_{i,j}$-s) such that $|{\hat Z}_n|>\gamma_2\alpha_2^n$ for all~$n$
with probability at least~$p_2$, where
$|{\hat Z}_n|={\hat Z}_n^0+\cdots+{\hat Z}_n^{t_0-1}$. In fact, we
think that with some more effort one should be able to prove that
these constants can be chosen uniformly in $q_i$-s and $q_{i,j}$-s;
however, it is easier to proceed as follows. Clearly, there are
$\gamma_3,p_3>0$, $\alpha_3>1$ ({\it not\/} depending on $q_i$-s and
$q_{i,j}$-s) such that
\begin{eqnarray}
\IP\Big[ \Po^0[|{\hat Z}_n|>\gamma_3\alpha_3^n] \geq p_3,\mbox{ for all }n
\;\Big\vert \;
0\in S_{\o}
\Big] &\geq & \frac{1}{2}.\phantom{***}
\label{ochen'_hor_seed}
\end{eqnarray}

Now, we recall the aperiodicity Condition~A.
Essentially, it says that the density of the aperiodic sites is positive,
where by ``aperiodic site'' we mean the following: for a given~$\o$,
$x\in\Z^d$ is an aperiodic site if there exists~$y$ such that
$\|x-y\|_1$ is even, and a particle in~$x$ sends at least one offspring to~$y$
with a positive $\Po$-probability.

For any $z\in S_{\o}$ and~${\hat Z}$ the process defined above
starting from $z$, define the event
\[
  E_z = \{|{\hat Z}_n|>\gamma_3\alpha_3^n , \mbox{ for all } n \}.
\]
Define
\begin{eqnarray*}
 M'_n &=& \{\forall y\in\K_{L_0 n\ln^{-1} n} \exists z\in S_{\o}:
\Po^z[E_z]\geq p_3, \|y-z\|_{\infty} \leq \alpha\ln n, \\
&&~~~\mbox{and there is an aperiodic site~$x_1$ such that }
            x_1+\A\subset\K_{L_0 n\ln^{-1} n}\}.
\end{eqnarray*}
Due to~(\ref{ochen'_hor_seed}), for $\IP[(M'_n)^c]$ we have an estimate
similar to~(\ref{term_2}). Analogously to the proof of
part~(i) of the lemma, one can prove that with overwhelming
probability before time $n\ln^{-1}n$ the random walk~$\xi_n$ will meet
a seed (centered, say, in~$z_0$) where an ``explosion'' (i.e., the event~$E_{z_0}$)
happens. Suppose that~$t_1$ is the moment when the Galton-Watson process
in~$z_0$ starts; we have $t_1\leq n\ln^{-1}n$.

Now, take an arbitrary $m\geq n$ and suppose that $\o\in M'_n$.
Supposing that~$n$ is large enough and~$\delta_1$ is small enough, there
exists~$k_0$ such that
\[
[t_1+k_0t_0, t_1+(k_0+1)t_0) \subset [m-4\delta_1 n, m-3\delta_1 n]
\]
and $k_0\geq\frac{m}{2t_0}$. Then, since the event~$E_{z_0}$ occurs,
there exists $t_2\in\{0,\ldots,t_0-1\}$ such that
${\hat Z}_{k_0}^{t_2} \geq \gamma_3\alpha_3^{k_0}/t_0$; i.e., there
are at least $\gamma_3\alpha_3^{k_0}/t_0$ ``unused particles'' (they
were not used in construction of the branching processes, so we do
not have any information about their future) at~$z_0$ at the moment
$t_1+k_0t_0+t_2$. Take any $x_0\in \K_{\delta_1 n}$ and suppose, for
definiteness, that $\|x_0-z_0\|_1$ is odd. Denote
$\hat t = m-(t_1+k_0t_0+t_2)$
(notice that $3\delta_1 n \leq \hat t \leq 4\delta_1 n$)
and consider two cases:

\smallskip
\noindent
{\sl Case 1:} $\hat t$ is odd.

\noindent
Then, by Condition~UE, any particle in~$z_0$ will send a descendant to~$x_0$
in time {\it exactly\/}~$\hat t$ with probability at least
$\eps_0^{\hat t}$.

\smallskip
\noindent
{\sl Case 2:} $\hat t$ is even.

\noindent
Here we will have to use the fact that on~$M'_n$ there exists an
aperiodic site somewhere in $\K_{L_0 n\ln^{-1} n}$. That is, when
going from~$z_0$ to~$x_0$ in time~$\hat t$, on the way we pass through the
aperiodic site, and this happens with probability at least
$C_6\eps_0^{\hat t-\ell}$.

\smallskip

So, in both cases we see that a particle in~$z_0$ will send a descendant to~$x_0$
in time~$\hat t$ with probability at least $C_7\eps_0^{\hat t}$ for some~$C_7>0$.
Recall that we dispose of at least $\gamma_3\alpha_3^{k_0}/t_0$
independent particles in~$z_0$, so the probability that at least one particle will
be in~$x_0$ at time~$m$ is at least
\begin{eqnarray*}
1-(1-C_7\eps_0^{\hat t})^{\gamma_3\alpha_3^{k_0}/t_0} &\geq &
      1-(1-C_7\eps_0^{4\delta_1 n})^{\gamma_3t_0^{-1}\alpha_3^{m/(2t_0)}}\\
&\geq & 1-\exp\Big\{-\frac{C_8\gamma_3}{t_0}
               \exp\Big[\frac{\ln\alpha_3}{2t_0}m -
                         4\delta_1 n\ln\eps_0^{-1}\Big]\Big\}.
\end{eqnarray*}
Choosing~$\delta_1$ small enough, it is straightforward to complete
the proof of the part~(ii).

\medskip
\noindent
{\bf Step 2.}
As for the part~(i), it can be proved analogously to the part~(ii),
by writing
\[
 \{T(0,x)\leq n\} \supset (\{T(0,x)=n\}\cup\{T(0,x)=n-1\})
\]
and noting that for handling of one of the events in the right-hand
side of the above display the Condition~A is unnecessary.
The proof of Lemma~\ref{hit_time} is completed.
\qed

Consider any $x_0\in\Z^d\setminus\{0\}$, and define a family
of random variables
\[
 Y^{x_0}(m,n) = T(mx_0,nx_0), \qquad 0\leq m < n.
\]
Let us point out that the sequence of
random variables $(Y^{x_0}(n-1,n),n=1,2,3,\ldots)$ is in general {\it not\/}
stationary (although they are of course identically distributed).
To see this, note that, conditioned on~$\o$, the random variables
$T(0,x_0)$ and $T(x_0,2x_0)$ are independent (because, recall the
construction of Section~\ref{sec:constr}, given $T(0,x_0)=r$,
the random variable $T(x_0,2x_0)$ depends on $v^{x,i}(n)$ for $n\geq r$), 
while $T(x_0,2x_0)$ and
$T(2x_0,3x_0)$ need not be so. Nevertheless, we will prove
that the above sequence satisfies the Strong Law of Large Numbers:
\begin{lmm}
\label{l_LLN}
Denote $\beta_{x_0}=\EA Y^{x_0}(0,1)$. Then for any $\eps>0$ there exists
$\theta_2=\theta_2(\eps)$ such that
\begin{equation}
\label{oc_LLN}
\PA\Big[\Big|\frac{1}{n}\sum_{i=1}^n Y^{x_0}(i-1,i)-\beta_{x_0}\Big|>\eps\Big]
            \leq \exp\{-\theta_2 \ln^d n\}.
\end{equation}
for all~$n$. In particular,
\begin{equation}
\label{eq_LLN}
 \frac{1}{n}\sum_{i=1}^n Y^{x_0}(i-1,i) \longrightarrow \beta_{x_0}
            \qquad \mbox{$\PA$-a.s.\ and in $L^p, p \geq 1$.}
\end{equation}
\end{lmm}

\noindent
{\it Proof.} Abbreviate $Y_i:=Y^{x_0}(i-1,i)$ and introduce the events
$G_i=\{Y_i<\sqrt{n}/(2L_0)\}$, $i=1,\ldots,n$. Suppose for simplicity
that~$\sqrt{n}$ is integer, the general case can be treated analogously.
Define the events
\begin{eqnarray*}
F &=& \Big\{\Big|\frac{1}{n}\sum_{i=1}^n Y_i-\beta_{x_0}\Big|>\eps\Big\},\\
F_i &=& \Big\{\Big|\frac{1}{\sqrt{n}}
             \sum_{j=1}^{\sqrt{n}} Y_{i+(j-1)\sqrt{n}}-\beta_{x_0}\Big|>\eps\Big\},
\end{eqnarray*}
$i=1,\ldots,\sqrt{n}$; since
\[
\frac{1}{n}\sum_{i=1}^n Y_i = \frac{1}{\sqrt{n}}\sum_{i=1}^{\sqrt{n}}
               \sum_{j=1}^{\sqrt{n}} Y_{i+(j-1)\sqrt{n}},
\]
we can write
\begin{equation}
\label{eq_ob'ed}
 F \subset \bigcup_{i=1}^{\sqrt{n}} F_i.
\end{equation}

Now, to bound from above the probability of a single event~$F_i$, we write
\begin{eqnarray*}
\PA[F_i] &\leq & \PA\Big[\Big|\frac{1}{\sqrt{n}}
             \sum_{j=1}^{\sqrt{n}} 
                Y_{i+(j-1)\sqrt{n}}{\bf 1}_{G_{i+(j-1)\sqrt{n}}}-
              \beta_{x_0}\Big|>\eps\Big]\\
  &&  {} + \PA[\mbox{there exists $j \leq \sqrt{n}$
such that $G_{i+(j-1)\sqrt{n}}^c$ occurs}]\\
 & =: & I_1 + I_2.
\end{eqnarray*}
By Theorem~\ref{t_upbound_ht}, with some $C_9>0$
\begin{equation}
\label{oc_T2}
I_2 \leq \sqrt{n}\PA[G_1^c] \leq \sqrt{n}\exp\{-C_9\ln^d  n\}.
\end{equation}

To bound the term~$I_1$, we note first that it is elementary
to obtain from Theorem~\ref{t_upbound_ht} that, for some $C_{10}$,
\begin{equation}
\label{oc_T1_1}
\beta_{x_0} - \EA Y_1 {\bf 1}_{G_1} \leq \exp\{-C_{10}\ln^d n\}.
\end{equation}
The key point here is that the random variables
$Y_{i+(j_1-1)\sqrt{n}}{\bf 1}_{G_{i+(j_1-1)\sqrt{n}}}$
and $Y_{i+(j_2-1)\sqrt{n}}{\bf 1}_{G_{i+(j_2-1)\sqrt{n}}}$ are
independent when $j_1\neq j_2$. Indeed, on the event
\[
\Big\{\max\{T((n_1-1)x_0,n_1x_0),T((n_2-1)x_0,n_2x_0)\}
        < \frac{|n_2-n_1|}{2L_0}\Big\}\;,
\]
the random variables  $T((n_1-1)x_0,n_1x_0)$ and
$T((n_2-1)x_0,n_2x_0)$ are functions of  $v^{x,i}(n)$-s
where the superscript $x$
belongs to nonintersecting subsets of~$\Z^d$).
Therefore, having in mind~(\ref{oc_T1_1})
and Theorem~\ref{t_upbound_ht},
to bound the term~$I_1$ from above we can use some Large Deviation
result for the sums of i.i.d.\ random variables without exponential
moments (see e.g.\ Corollary~1.11 from~\cite{Nagaev}) to obtain that
\begin{equation}
\label{oc_T1}
I_1 < \exp\{-C_{11}\ln^d n\}.
\end{equation}
Using~(\ref{oc_T2}), (\ref{oc_T1}) and~(\ref{eq_ob'ed}), we conclude the
proof of~(\ref{oc_LLN}). Since~(\ref{eq_LLN}) follows from~(\ref{oc_LLN})
immediately for $p=1$, the proof of Lemma~\ref{l_LLN} is finished in
this case. To extend it to a general~$p$, it suffices to note
that for all $p'\geq 1$,
\[
\Big(\frac{1}{n}\sum_{i=1}^n Y^{x_0}(i-1,i)\Big)^{p'} \leq
\frac{1}{n}\sum_{i=1}^n (Y^{x_0}(i-1,i))^{p'},
\]
which has a finite expectation.
\qed

To proceed with the proof of Theorem~\ref{shape_t},
we state the result of~\cite{L},
 which is an improved version of Kingman's subadditive ergodic
theorem~\cite{K}.

\begin{theo}
\label{t_subadd}
Suppose that $\{Y(m,n)\}$ is a collection of
positive random variables indexed by integers satisfying $0\le m < n $
such that
\begin{itemize}
\item[(i)] $ Y(0,n) \le Y(0,m) + Y(m,n)$ for all $ 0 \le m < n $ (subadditivity);
\item[(ii)] the joint distribution of $\{Y(m+1, m+k+1), k \ge 1\}$
is the same as that of $\{ Y(m, m+k), k \ge 1\}$ for each $m\ge 0$;
\item[(iii)] for each $k \ge 1$ the sequence of random variables
$\{Y(nk,(n+1)k), n\ge 0\}$ is a stationary ergodic process;
\item[(iv)] The expectation of $Y(0,1)$ is finite.
\end{itemize}
Then
\[
 \lim_{n \to \infty} \frac{Y(0,n)}{n} \to \gamma \qquad \mbox{a.s.},
\]
where
\[
 \gamma = \inf_{n \ge 0} \frac{\EA Y(0,n)}{n}.
\]
\end{theo}

Similarly to the proof of a number of other shape results, our original
intention was to apply Theorem~\ref{t_subadd} to the family
$(Y^{x_0}(m,n), 0\leq m<n)$. Indeed, the assumption~(i) of
Theorem~\ref{t_subadd} holds due to Lemma~\ref{l_subadd},
  from the construction of the random variables $T(\cdot,\cdot)$
is is elementary to observe that the assumption~(ii) holds
as well, and the assumption~(iv) follows from Theorem~\ref{t_upbound_ht}.
However, as we observed just before Lemma~\ref{l_LLN}, the sequence
of random variables in~(iii) need not be stationary (even though it
has good mixing properties and the random variables there are
equally distributed). So, we take a slightly different route:
consider the proof
of Theorem~\ref{t_subadd} (here we use the proof of Theorem~2.6 of
Chapter~VI of~\cite{L_book}), and follow its steps carefully.
One sees that
the assumption~(iii) (which is the assumption~(b) in Theorem~2.6 of
Chapter~VI of~\cite{L_book})
is used only in~(2.11) and between the displays~(2.14)
and~(2.15) of Chapter~VI of~\cite{L_book} to prove that
a certain sequence converges a.s.\ and in~$L^1$ to its mean; in
our situation, Lemma~\ref{l_LLN} takes care of that.

From the above argument we conclude that for any $x\in\Z^d\setminus\{0\}$
there exists a number $\mu(x)$ (depending also on~$Q$) such that
\begin{equation}
\label{principal_shape}
\frac{T(0,nx)}{n} \longrightarrow \mu(x) \qquad \mbox{$\PA$-a.s., $n\to\infty$.}
\end{equation}

   From this point on, the proof of the shape result for~$B_n^0$ becomes
completely standard, so we only briefly outline the main steps and refer
to e.g.~\cite{AMP,BG,DG} for details:
\begin{itemize}
\item it is easy to obtain that for any $x\in\Z^d, a\in\Z_+$,
we have $\mu(ax)=a\mu(x)$;
\item using that, $\mu(x)$ is first extended on $x\in\R^d$ with rational
coordinates (if $x\in\Q^d$ and $ax\in\Z^d$, with $a\in\Z_+$,
then $\mu(x):=\frac{\mu(ax)}{a}$),
and then, using the subadditivity and the part~(i) of Lemma~\ref{hit_time},
to the whole~$\R^d$;
\item the limiting shape~$B$ is then identified by $B=\{x\in\R^d:\mu(x)\leq 1\}$
(note that~$B$ is convex since the subadditivity 
property $\mu(x+y)\leq \mu(x)+\mu(y)$
is preserved; however, $B$ need not be symmetric, since generally
$\mu(x)$ need not be equal to $\mu(-x)$);
\item to complete the proof (for $B_n^0$), 
cover~$B$ and a sufficiently large annulus
of~$B$ by balls of radius~$\delta'$, where $\delta'$ is sufficiently small,
and then use the part~(i) of Lemma~\ref{hit_time}.
\end{itemize}

To complete the proof of Theorem~\ref{shape_t} 
(for dimension $d\geq 2$), we recall the relation
${\tilde B}^x_n \subset {\bar B}^x_n  \subset B^x_n$, so all we need to
prove is that for any~$\eps>0$, $(1-\eps)B\subset\FF({\tilde B}^0_n)$
for all~$n$ large enough. This follows easily from
the part~(ii) of Lemma~\ref{hit_time} and the corresponding shape result
for~$B^0_n$.
\qed

\subsection{Case $d=1$}
\label{s_d=1}
As noticed in the beginning of Section~\ref{pf-shape_th},
here we cannot guarantee that $\EA T(0,1)<\infty$ 
(although it may be so), so we need to develop a 
different approach. On the other hand, still a number of
the steps of the proof for~$d=1$ will be quite analogous to
the corresponding steps of the proof for $d\geq 2$; in such cases
we will prefer to refer to the case $d\geq 2$ rather than
writing down a similar argument once again. 

The main idea of the proof of Theorem~\ref{shape_t}
in the case $d=1$ is the following. From the proofs
of Theorems~\ref{t_upbound_ht} and~\ref{t_lobound_ht} we saw
that, while usually $\Po[T(0,1)>n]$ is well behaved
(and, in particular, $\Eo T(0,1)<\infty$), there are some
``exceptional'' environments that may cause $\EA T(0,1)=\infty$
in dimension~$1$ (see Example~\ref{ex_infexp}).
So, if the environment is ``untypical'', instead of
starting with one particle, we start with a number of
particles depending on the environment (and the more untypical
is the environment, the larger is that number). 

For the sake of simplicity, we suppose now that the maximal jump~$L_0$
is equal to~$1$; afterwords we explain how to deal with general~$L_0$.

Keeping the notation~$S_\o$ from the proof of Theorem~\ref{t_upbound_ht},
we note that the set~$S_\o$ has positive density in~$\Z$, so
there exists (small enough) $\gamma_1$ such that an interval
of length~$k$ contains at least~$\gamma_1 k$ good seeds from~$S_\o$
with $\IP$-probability at least $1-e^{-C_1 k}$. Let us say that
an interval is {\it nice}, if (being~$k$ its length) it
contains at least~$\gamma_1 k$ good seeds from~$S_\o$.

Fix $r<C_1$ (e.g., $r:=C_1/2$) and define
\begin{eqnarray*}
 h^x(\o) &=& \min\{\text{$m$: all the intervals of length $k\geq m$} \\
 && ~~~~~~~~~~~~~~~~~~~~~~
 \text{intersecting with $x+[-e^{rk},e^{rk}]$ are nice}\},
\end{eqnarray*}
with this choice of~$r$ it is elementary to obtain that
there exists~$C_2>0$ such that
\begin{equation}
\label{hvost_odnako}
 \IP[h^0(\o) = n] \leq e^{-C_2 n}.
\end{equation}

Now, suppose that, instead of starting with one particle, the process
starts with $e^{Kh^0(\o)}$ particles in~$0$, where~$K$ is a (large)
constant to be chosen later. For $\ell\geq 1$, define
\[
 \tT(0,\ell) = \min\{n\geq 0 : \eta^0_n(\ell) \geq e^{Kh^\ell(\o)}\},
\]
i.e., $\tT(0,\ell)$ is the first moment when we have at least
$e^{Kh^\ell(\o)}$ particles in~$\ell$. Now, our goal is to prove
that if~$K$ is large enough, then $\EA\tT(0,1)<\infty$.
Denote $Z=h^0(\o)\vee h^1(\o)$, and write
\begin{eqnarray}
\EA\tT(0,1) &=& \IE \Eo\tT(0,1) \nonumber\\
 &\leq &\sum_{m=1}^\infty \Big(\sup_{\o:Z=m}\Eo\tT(0,1)\Big) 
                                                      \IP[Z=m]. 
\label{conta_d=1}
\end{eqnarray}

Let us obtain an upper bound on the supremum in the right-hand side
of~(\ref{conta_d=1}). Fix~$m\geq 1$ and let us consider an
environment~$\o$ such that $Z=m$ (so that $h^0(\o)\leq m$). 
First, we prove the estimate~(\ref{hvost_fiks_omega}) below,
in the following way:
\begin{itemize}
\item[(i)] Consider the time interval $[0,\theta_0 m]$, 
where $\theta_0=\frac{K}{2\ln\eps_0^{-1}}$. 
Each particle that is initially in the origin
(even if it does not generate new offsprings)
will cover the box $[0,\theta_0 m]\subset\Z$ by time $\theta_0 m$
(simply by going always one unit to the right), with probability
at least~$\eps_0^{\theta_0 m}$ ($\eps_0$ is from Condition~UE).
Recall that initially we had $e^{Kh^0(\o)}$ particles in~$0$;
since $K>\theta_0\ln\eps_0^{-1}$, there exists~$C_3$
such that, with probability at least $1-e^{-C_3m}$ all the sites of the
box $[0,\theta_0 m]$ will be visited by time $\theta_0 m$.
\item[(ii)] By definition of the quantity $h^0(\o)$, the
box $[0,\theta_0 m]$ contains at least $\theta_0 \gamma_1 m$
good seeds from~$S_\o$ (here we suppose that $\theta_0>1$, 
i.e., $K>2\ln\eps_0^{-1}$). Since all of them were visited, 
there will be an explosion in at least one of these seeds
with probability at least $1-e^{-C_4\theta_0m}$. 
\item[(iii)] Now, we only have to wait $C_5 m$ (where $C_5$ is a
[large] constant depending on~$K$) time units more
to be able to guarantee that at least $e^{Kh^1(\o)}$ particles
will simultaneously be in the site~$1$ at some moment from
the time interval $[\theta_0 m, \theta_0m+C_5m]$ (to see that,
use an argument of the type ``if  the number of visits to the site~$1$
during the time interval of length~$n_1$ was at least~$n_2$,
then at some moment at least $n_2/n_1$ particles were simultaneously in
that site'').
So, finally one can obtain that there exist $C_6, \theta_1$
(depending on~$K$) such that
\begin{equation}
\label{hvost_fiks_omega}
 \Po[\tT(0,1)>C_6 m] < e^{-\theta_1 m},
\end{equation}
and the crucial point is that~$\theta_1$ can be made
arbitrarily large by enlarging~$K$. So, choose~$K$ in such a way that
$\theta_1>5\ln\eps_0^{-1}$.
\end{itemize}

Next, the goal is to obtain an upper estimate on $\Po[\tT(0,1)>n]$
which does not depend on~$K$. Specifically, we are going to 
prove that, for some positive constants $C_7,C_8$, we have,
on $\{\o : h^0(\o)\vee h^1(\o)=m\}$ and for large enough~$m$
\begin{equation}
\label{oc_strexp}
 \Po[\tT(0,1)>n] \leq e^{-C_7 n^{C_8}}
\end{equation}
for all $n\geq e^{5m\ln\eps_0^{-1}}$. 

\medskip
\noindent
{\bf Remark.}
To obtain the estimate~(\ref{oc_strexp}),
we will use only one initial particle in~$0$; so, the same estimate will
be valid for~$T(0,1)$, thus giving us the proof of Proposition~\ref{quench_dim1}.
\qed

\medskip

Now, to prove~(\ref{oc_strexp}), we proceed in the following way.
\begin{itemize}
\item[(i)] Consider one particle starting from the origin.
During any time interval of length $\frac{\ln n}{5\ln\eps_0^{-1}}$
it will cover a space interval of the same length (by going to the right
on each step) with probability at least
\[
\eps_0^{\frac{\ln n}{5\ln\eps_0^{-1}}}=n^{-1/5}
\]
(note that there is a similar argument in the proof of 
Theorem~\ref{t_upbound_ht} for the case $d=1$).
So, in time~$n^{1/4}$ a single particle will cover an
interval of that length with probability at least $1-e^{-C_9n^{1/20}}$
(note that these estimates do not depend on~$\o$).
\item[(ii)] Abbreviate $r'=\frac{r}{5\ln\eps_0^{-1}}$
($r$ is from the definition of~$h^x(\o)$). If $n\geq e^{5m\ln\eps_0^{-1}}$,
then all the intervals of length $\frac{\ln n}{5\ln\eps_0^{-1}}$ intersecting
with the interval $[-n^{r'},n^{r'}]$ are nice (so, in particular, they
contain at least one good seed from~$S_\o$)
on $\{\o : h^0(\o)\vee h^1(\o)=m\}$. 
\item[(iii)] Consider the time interval $[0,n^{1/2}]$. One of the 
following two alternatives will happen:
 \begin{itemize}
  \item[(iii.a)] Either some of the particles from the 
             cloud of the offsprings of the initial 
             particle will go out of the interval $[-n^{r'},n^{r'}]$, or
  \item[(iii.b)] all the offsprings of the initial particle
             will stay in the interval $[-n^{r'},n^{r'}]$ up to time~$n^{1/2}$.
 \end{itemize}
In the case (iii.a), at least $\gamma_1 n^{r'}$ good seeds from~$S_\o$
will be visited. In the case (iii.b), argue as follows: we have $n^{1/4}$
time subintervals of length~$n^{1/4}$; during each one a good seed will
be visited with overwhelming probability. So, with probability
greater than $1-n^{1/4}e^{-C_9n^{1/20}}$
the number of visits to good seeds will be at least $n^{1/4}$
(and all of these good seeds are in the interval $[-n^{1/2},n^{1/2}]$).
\item[(iv)] Thus, in any case, by time $n^{1/2}$ there will be
a polynomial number of visits to good seeds. So, with overwhelming probability
one of them will explode and produce enough particles to guarantee
that there are at least $e^{C_{10}n}$ particles which were
created at distance no more than~$n^{1/2}$ from~$0$ before time $n/2$.
Then, it is elementary to obtain that, with overwhelming probability,
we will have at least $e^{C_{11}n}$ particles in site~$1$, for some $C_{11}>0$.
Since $n\geq e^{5m\ln\eps_0^{-1}}$, this will be enough to
to make the event $\{\tT(0,1)\leq n\}$ occur when $C_{11}e^{5m\ln\eps_0^{-1}}>Km$,
and the last inequality holds for all, except possibly finitely many,
values of~$m$.
\end{itemize}

Now, we finish the proof of the fact that $\EA\tT(0,1)<\infty$.
Write
\[
 \Eo \tT(0,1) = \sum_{n=0}^\infty \Po[\tT(0,1) > n],
\]
and use~(\ref{oc_strexp}) to bound $\Po[\tT(0,1) > n]$ for
$n\geq e^{5m\ln\eps_0^{-1}}$, and~(\ref{hvost_fiks_omega})
for $n\in [C_6 m, e^{5m\ln\eps_0^{-1}})$. 
Since $\theta_1>5\ln\eps_0^{-1}$, we obtain that
$\Eo \tT(0,1) < C_{12}m + C_{13}$ 
on $\{\o : h^0(\o)\vee h^1(\o)=m\}$ for~$m$ large enough, so
using~(\ref{conta_d=1}) and~(\ref{hvost_odnako}), we
conclude the proof of the fact that $\EA\tT(0,1)<\infty$.

The rest of the proof (for $L_0=1$) is straightforward.
First, we define variables $\tT(k,m)$, $1\leq k<m$,
repeating the construction of Section~\ref{sec:constr},
with the following modifications: the process initiating in~$k$
starts from $e^{Kh^k(\o)}$ particles, at the moment (with
respect to $v^{x,i}(n)$) $\tT(0,k)$ (instead of $T(0,k)$).
Then, it is elementary to see that we still have the subadditivity relation
$\tT(0,m)\leq \tT(0,k)+\tT(k,m)$. There is again a problem with the
absence of the stationarity for the sequence $\tT(0,1), \tT(1,2), \tT(2,3)\ldots$;
this problem can be dealt with in exactly the 
same way as in Section~\ref{s_d_geq2}.

So, the above arguments show that $\frac{\tT(0,n)}{n}$ converges to
a limit as $n\to\infty$, which immediately implies the shape theorem
in dimension~$1$ (we do not need the analogue of Lemma~\ref{hit_time}
here) for the model starting with $e^{Kh^0(\o)}$ particles
from~$0$. 

We now complete the proof of Theorem~\ref{shape_t} in the case $d=1$
(and, for now, $L_0=1$): it is elementary to
obtain that, for a recurrent branching random walk in random environment
starting with~$1$ particle, for $\IP$-almost all $\o$-s, at some
(random) time we will have at least $e^{Kh^0(\o)}$ particles in the origin.
Now, it remains only to erase all other particles and apply the above
reasoning.

To treat the case of a general $L_0\geq 1$, we apply the same
reasoning as in the proof of Theorem~\ref{t_upbound_ht} for $d=1$
(namely, instead of $(\eps,\rho)$-good seeds, we consider 
$(\eps,\rho,W)$-good seeds with $W=\{0,1,\ldots,L_0-1\}$, so
that a particle cannot overjump~$W$).
\qed

\end{document}